\documentclass[12pt,leqno]{amsart}

\usepackage{float}
\usepackage{amssymb,amsmath,amsthm,amscd}
\usepackage{graphics}
\usepackage{mathrsfs}
\usepackage{enumerate}
\usepackage{color}
\usepackage{verbatim}
\usepackage[bookmarks=true,bookmarksnumbered=true,bookmarkstype=toc=true,pdftitle={\@title},pdfauthor={\@author},pdfstartview={FitBH -32768}]{hyperref}

\usepackage[all]{xy}
\numberwithin{equation}{section}
\newcommand{\nc}{\newcommand}
\usepackage{graphicx}
\usepackage{tikz}

\def\node#1#2{\overset{#1}{\underset{#2}{\circ}}}

\def\ver#1#2{\overset{{\llap{$\scriptstyle#1$}\displaystyle\circ{\rlap{$\scriptstyle#2$}}}}{\scriptstyle\vert}}

\usetikzlibrary{arrows}
\tikzstyle{every picture}+=[remember picture]
\tikzstyle{na} = [baseline=-.5ex]
\tikzstyle{mine}= [arrows={angle 90}-{angle 90},thick]

\def\Llleftarrow{%
\lower2pt\hbox{\begingroup
\tikz
\draw[shorten >=0pt,shorten <=0pt] (0,3pt) -- ++(-1em,0) (0,1pt) -- ++(-1em-1pt,0) (0,-1pt) -- ++(-1em-1pt,0) (0,-3pt) -- ++(-1em,0) (-1em+1pt,5pt) to[out=-105,in=45] (-1em-2pt,0) to[out=-45,in=105] (-1em+1pt,-5pt);
\endgroup}
}

\def\Rrrightarrow{\vcenter{\hbox{\rotatebox{180}{\Llleftarrow}}}}

\newcommand{\BR}{\boldsymbol{r}}
\newcommand{\KORP}{\mathcal{P}}
\newcommand{\KORC}{\mathcal{C}}
\newcommand{\CART}{\mathsf{Cart}}
\newcommand{\KORR}{\mathcal{R}}

\newcommand{\NNN}{\N_{+}}
\newcommand{\PPP}{\Psi}
\newcommand{\emptypartition}{\phi}
\newcommand{\SORT}{\mathsf{Sort}}
\newcommand{\IRR}{\mathsf{Irr}}
\newcommand{\SUPER}{\mathsf{su}}
\newcommand{\HT}{\mathsf{ht}}

\newcommand{\PAR}{\mathsf{Par}}
\newcommand{\RPAR}{\mathsf{RP}}
\newcommand{\CPAR}{\mathsf{CRP}}

\newcommand{\WI}{\boldsymbol{i}}

\newcommand{\ACT}{\odot}
\newcommand{\XXX}{\boldsymbol{x}}

\newcommand{\MMP}{\widetilde{P}^{-}}

\newcommand{\C}{{\mathbb C}}
\newcommand{\BL}{S}
\newcommand{\BLOCK}{\mathsf{Bl}}
\newcommand{\BAR}{{\tau}}
\newcommand{\GR}{\mathsf{Gd}}
\newcommand{\GRM}{\mathsf{GM}}
\newcommand{\SH}{\mathsf{Sh}}
\newcommand{\SIM}{\sim}
\newcommand{\QSH}{\mathsf{QSh}}
\newcommand{\RSH}{\mathsf{RSh}}
\newcommand{\SHM}{\mathsf{Sh}^\mathsf{M}}
\newcommand{\QSHM}{\mathsf{QSh}^\mathsf{M}}
\newcommand{\RSHM}{\mathsf{RSh}^\mathsf{M}}
\newcommand{\Q}{\mathbb {Q}}

\newcommand{\Z}{{\mathbb Z}}

\newcommand{\MCC}{{\mathscr{C}}}
\newcommand{\MP}{{\mathcal{P}}}
\newcommand{\MPC}{{\mathcal{P}^\vee}}
\newcommand{\MA}{\mathscr{A}}
\newcommand{\MAAA}{\mathscr{H}}
\newcommand{\MAA}{{R}}
\newcommand{\corr}{\mathbb{F}}

\newcommand{\TRANS}[1]{{}^{\textrm{tr}}{#1}}
\newcommand{\corrr}[1]{k_{#1}}
\newcommand{\qq}[1]{\eta_{#1}}

\theoremstyle{plain}
\newtheorem{lemma}{Lemma}[section]
\newtheorem{prop}[lemma]{Proposition}
\newtheorem{theorem}[lemma]{Theorem}
\newcommand{\Prop}{\begin{prop}}
\newcommand{\enprop}{\end{prop}}
\newcommand{\Lemma}{\begin{lemma}}
\newcommand{\enlemma}{\end{lemma}}
\newcommand{\Th}{\begin{theorem}}
\newcommand{\enth}{\end{theorem}}
\newtheorem{corollary}[lemma]{Corollary}
\newcommand{\Cor}{\begin{corollary}}
\newcommand{\encor}{\end{corollary}}
\newtheorem{definition}[lemma]{Definition}
\newtheorem{question}[lemma]{Question}
\newcommand{\Question}{\begin{question}}
\newcommand{\enquestion}{\end{question}}
\newcommand{\Def}{\begin{definition}}
\newcommand{\edf}{\end{definition}}

\theoremstyle{definition}
\newtheorem{remark}[lemma]{Remark}

\newtheorem{conjecture}[lemma]{Conjecture}
\nc{\Rem}{\begin{remark}}
\nc{\enrem}{\end{remark}}

\newcommand{\Con}{\begin{conjecture}}
\newcommand{\encon}{\end{conjecture}}

\nc{\emprule}[1]{\rule{#1}{0pt}}
\newcommand{\g}{{\mathfrak{g}}}
\newcommand{\Gg}{{\mathfrak{g}}}
\newcommand{\Ggz}{{\mathfrak{g}_0}}
\newcommand{\Ggg}{{\widehat{\mathfrak{g}}_0}}

\newcommand{\Hom}{\operatorname{Hom}}
\newcommand{\Rank}{\operatorname{rank}}

\newcommand{\isoto}[1][]{\mathop{\xrightarrow[#1]{\rule{0pt}{.9ex}{\raisebox{-.6ex}[0ex][-.7ex]{$\mspace{4mu}\sim\mspace{3mu}$}}}}}

\renewcommand{\hom}{\operatorname{\it \mathscr{H}\kern-.25em om}}

\newcommand{\KSS}[2]{{\{{#1}\}_{{#2}}}}
\newcommand{\KSSS}[3]{{\{{#1}\}^{{#3}}_{{#2}}}}
\newcommand{\Binom}[3]{{{#1} \atopwithdelims[] {#2}}_{#3}}
\newcommand{\Bino}[2]{{{#1} \atopwithdelims[] {#2}}}

\newcommand{\M}{{\mathscr M}}
\newcommand{\N}{{\mathbb{N}}}
\newcommand{\eq}{\begin{eqnarray}}
\newcommand{\eneq}{\end{eqnarray}}
\newcommand{\GPPP}[1]{\mathsf{Proj}_\mathsf{gr}({#1})}
\newcommand{\GP}[1]{\mathsf{Proj}({#1})}

\newcommand{\eqn}{\begin{eqnarray*}}
\newcommand{\eneqn}{\end{eqnarray*}}

\newenvironment{tenumerate}{
  \begin{enumerate}
  
  }{\end{enumerate}}

\nc{\bnum}{\begin{enumerate}[{\rm (i)}]}
\nc{\enum}{\ee}
\nc{\bna}{\begin{enumerate}[{\rm (a)}]}
\nc{\ena}{\ee}

\newcommand{\on}{\operatorname}
\newcommand{\bni}{\begin{tenumerate}}
\newcommand{\eni}{\end{tenumerate}}

\newcommand{\QED}{\end{proof}}
\newcommand{\Proof}{\begin{proof}}

\newcommand{\CONG}[1]{\equiv_{#1}}

\newcommand{\ba}{\begin{array}}
\newcommand{\ea}{\end{array}}

\newcommand{\hs}{\hspace*}

\newcommand{\eqsub}{\begin{subequations}\begin{eqnarray}}
\newcommand{\eneqsub}{\end{eqnarray}\end{subequations}}

\newcommand{\ol}{\overline}

\newcommand{\A}{\mathscr{A}}

\nc{\la}{\lambda}
\nc{\lam}{\lambda}
\nc{\U}[1][\g]{U_q(#1)}
\nc{\te}{\tilde{e}}
\nc{\tei}{\tilde{e}_i}
\nc{\tf}{\tilde{f}}
\nc{\tfi}{\tilde{f}_i}
\nc{\tU}{\widetilde U_q(\g)}
\nc{\tE}{\tilde{E}}
\nc{\tF}{\tilde{F}}

\nc{\BZ}{{\mathbb{Z}}}
\nc{\al}{\alpha}
\nc{\qs}{{q}}
\nc{\lan}{\langle}
\nc{\ran}{\rangle}
\nc{\re}{{\mathrm{re}}}
\nc{\wt}{\operatorname{wt}}
\nc{\Uf}[1][\g]{U^-_q(#1)}
\nc{\Ue}{U^+_q(\g)}
\nc{\eps}{\varepsilon}
\nc{\vphi}{\varphi}
\nc{\sphi}{\varphi^*}
\nc{\seps}{\varepsilon^*}

\nc{\nn}{\nonumber}

\nc{\vp}{\varpi}
\nc{\cls}{{\operatorname{cl}}}

\nc{\Wt}{{\operatorname{Wt}}}
\nc{\Us}{U'_q(\g)}
\nc{\La}{\Lambda}
\nc{\ro}{{\rm(}}
\nc{\rf}{{\rm)}}
\nc{\norm}{{\mathrm{norm}}}
\nc{\qbox}{\quad\mbox}
\nc{\braid}{{\mathfrak{B}}}
\nc{\Ad}{\operatorname{Ad}}
\nc{\Aut}{\operatorname{Aut}}
\nc{\dt}[1]{\tilde{\tilde #1}}
\nc{\Sn}{S^{{\mathrm{norm}}}}
\nc{\aff}{{\mathrm{aff}}}
\nc{\rk}{{\mathrm{rk}}}
\nc{\tQ}{\widetilde{Q}}
\nc{\tP}{\widetilde{P}}
\nc{\tW}{\widetilde{\mathscr{W}}}
\nc{\Dyn}{\mathrm{Dyn}}
\nc{\tD}{\widetilde{\Delta}}
\nc{\height}{{\operatorname{ht}}}
\nc{\bl}{\bigl}
\nc{\br}{\bigr}
\nc{\Hecke}{\mathrm{H}}
\nc{\HB}{\Hecke^{\mathrm{B}}}
\nc{\K}{\mathrm{K}}
\newcommand{\scbul}{{\,\raise1pt\hbox{$\scriptscriptstyle\bullet$}\,}}
\nc{\vac}{{\phi}}
\nc{\be}{\begin{enumerate}}
\nc{\ee}{\end{enumerate}}
\nc{\low}{{\mathrm{low}}}
\nc{\upper}{{\mathrm{up}}}
\nc{\Zodd}{\Z_{\mathrm{odd}}}
\nc{\Ft}[1][n]{\mathbb{P}\mathrm{ol}_{#1}}
\nc{\Ftf}[1][n]{\widetilde{\mathbb{P}\mathrm{ol}}_{#1}}
\nc{\KA}{\on{K}^{\mathrm{A}}}
\nc{\KB}{\on{K}^{\mathrm{B}}}
\nc{\Res}{\on{Res}}
\nc{\Fc}[1][{n,m}]{\mathbf{F}_{#1}}
\nc{\tphi}{\widetilde{\varphi}}
\nc{\CO}{\mathscr{O}}
\nc{\CK}{\mathscr{K}}
\nc{\disc}{\mathfrak{d}}
\nc{\tr}{\on{tr}}
\nc{\Gb}{\mathfrak{b}}
\nc{\Gh}{\mathfrak{h}}
\nc{\ga}{\mathfrak{a}}
\nc{\stable}{\mathrm{stable}}
\nc{\X}{\mathfrak{X}}
\nc{\Hilb}{\mathrm{Hilb}}
\nc{\W}{\ensuremath{\mathscr{W}}}
\nc{\Ws}{\ensuremath{\rm W}}
\nc{\opp}{{\on{opp}}}

\nc{\corps}{{\mathbf{k}}}
\nc{\cor}{{\mathbf{k}}}
\nc{\h}{\mathrm{\hslash}}
\nc{\fL}[1][{\h}]{\C(\mspace{-1mu}(#1)\mspace{-1mu})}
\nc{\ad}{\mathrm{ad}}
\newcommand{\Endm}{\operatorname{\mathscr{E}\kern-.1pc\mathit{nd}}}
\newcommand{\Endomo}{\operatorname{\mathscr{E}\kern-.1pc\mathit{nd}}}
\nc{\bc}{\bar{\corps}}
\nc{\reg}{{\mathrm{reg}}}
\nc{\ysq}{\mathbf{y}^2}
\nc{\CH}{\mathsf{char}}

\nc{\sketch}{\Proof}
\nc{\Gm}{\mathbb{G}_{\mathrm{m}}}
\nc{\hGm}{\hat{\mathbb{G}}_{\mathrm{m}}}
\nc{\ug}{\widehat{\mathrm{G}}_{\mathrm{m}}}

\nc{\tL}{\widetilde{\mathscr{L}}}
\nc{\Fr}{\mathcal{F}}
\nc{\E}{\mathcal{E}}

\nc{\ord}{\on{ord}}
\nc{\bM}{\overset{\hs{1.5ex}\rule[-.08ex]{1.8ex}{.08ex}}{\M}}
\nc{\romano}{\mathrm{o}}
\nc{\into}{\hookrightarrow}
\nc{\good}{\mathrm{good}}
\nc{\tA}{\widetilde\A}
\nc{\Vz}{{V}\kern-1.1ex\raisebox{1.5ex}[0ex][0ex]{$\cdot$}}
\nc{\bxes}[1]{\raisebox{.9ex}{$\cdot$}{\kern#1}\raisebox{0ex}{$\cdot$}
{\kern#1}\raisebox{-.9ex}{$\cdot$}}
\nc{\ssum}{\mathop{\mbox{\normalsize{${\sum}$}}}\limits}
\nc{\ct}{{\mbox{\tiny$\mathrm{CT}$}}}
\nc{\pr}{\mathrm{pr}}
\nc{\qr}{\mathrm{rp}}
\DeclareMathOperator{\SPAN}{\mathsf{span}}
\DeclareMathOperator{\RAD}{\mathsf{Rad}}
\DeclareMathOperator{\MAT}{\mathsf{Mat}}
\DeclareMathOperator{\END}{\mathsf{End}}

\DeclareMathOperator{\COKER}{\mathsf{Coker}}

\DeclareMathOperator{\SMT}{\mathsf{SIM}}
\DeclareMathOperator{\SYM}{\mathsf{Sym}}
\DeclareMathOperator{\ID}{\mathsf{id}}
\DeclareMathOperator{\ORD}{\mathsf{ord}}
\DeclareMathOperator{\GL}{\mathsf{GL}}
\DeclareMathOperator{\SYL}{\mathsf{Syl}}
\DeclareMathOperator{\PRS}{\mathsf{pdiv}}
\DeclareMathOperator{\GCD}{\mathsf{gcd}}
\DeclareMathOperator{\PROC}{\mathsf{PC}}
\DeclareMathOperator{\DIAG}{\mathsf{diag}}
\DeclareMathOperator{\RED}{\mathsf{red}}
\DeclareMathOperator{\RES}{\mathsf{res}}
\DeclareMathOperator{\CUT}{\mathsf{cut}}
\DeclareMathOperator{\INFL}{\mathsf{Infl}}

\def\iso{\isoto}

\def\GL{\operatorname{GL}\nolimits}

\nc{\Fs}{\ensuremath{\rm F}}

\nc{\isotf}{\overset{
{\rule{0pt}{.9ex}%
{\raisebox{-.6ex}[0ex][-.7ex]{$\mspace{3mu}\sim\mspace{3mu}$}}}}
{\longleftrightarrow}}
\nc{\tN}{\tilde\N}
\nc{\tens}{\mathop\otimes\limits}
\nc{\super}{\mathrm{super}}
\nc{\Mods}{\on{Mod}_\super}
\nc{\rev}{\mathrm{rev}}
\nc{\clif}{\mathfrak{C}}
\nc{\clifm}{\clif^{-}}
\nc{\Fct}{\mathrm{Fct}}
\nc{\Fcts}{\mathrm{Fct}_\super}

\nc{\Ks}{\on{K^\super}}
\nc{\ts}{\widetilde{s}}
\nc{\KUGIRI}{\circ}
\nc{\Sym}{\mathfrak{S}}
\nc{\FF}{\mathcal{F}}
\nc{\BF}{\mathcal{B}}
\nc{\BFF}{\widetilde{\mathcal{B}}}

\nc{\cc}{\mathfrak{c}}
\nc{\SK}{\mathcal{KS}}
\nc{\noi}{\noindent}
\nc{\odd}{{\mathrm{odd}}}
\nc{\even}{{\mathrm{even}}}
\nc{\bs}{\ol{s}}
\nc{\Khc}[1][n]{\ol{\mathcal{KHC}}_{#1}}
\nc{\Ohc}[1][n]{\ol{\mathcal{OHC}}_{#1}}
\nc{\KHC}[1][n]{\mathcal{K}{\mathcal{HC}}_{#1}}
\nc{\OHC}[1][n]{\mathcal{O}{\mathcal{HC}}_{#1}}
\nc{\IODD}{I_{\odd}}
\nc{\IEVEN}{I_{\even}}
\DeclareMathOperator{\KKK}{\mathsf{K}_0}
\newcommand{\MOD}[1]{{\mathsf{Mod}({#1})}}
\newcommand{\PROJ}[1]{{\mathsf{Proj}({#1})}}
\newcommand{\GMOD}[1]{{\mathsf{Mod}_{\mathsf{gr}}({#1})}}
\newcommand{\F}{\mathbb{F}}
\nc{\IRED}{I}
\DeclareMathOperator{\CHAR}{char}

\nc{\MH}{\mathcal{H}}

\newcommand{\KLR}{R}

\newcommand{\MAPSTO}{\longmapsto}

\newcommand{\SEQ}{\mathsf{Seq}}
\DeclareMathOperator{\DEG}{\mathsf{deg}}

\newcommand{\HA}{\widehat{X}}
\newcommand{\HI}{\widehat{I}}

\newcommand{\RP}[1]{\mathsf{RP}_3}

\nc{\bwr}{\mbox{\large$\wr$}}
\nc{\At}[1][{{i,j}}]{\mathscr{A}_{{#1}}}
\nc{\tAt}[1][{{i,j}}]{\widetilde{\mathscr{A}}_{{#1}}}
\nc{\Bt}[1][{{i,j}}]{\mathscr{B}_{{#1}}}
\nc{\tBt}[1][{{i,j}}]{\widetilde{\mathscr{B}}_{{#1}}}
\nc{\prt}[1]{\mathrm{par}(#1)}
\nc{\er}{\mathrm{e}}
\nc{\ec}{\mathrm{e}^-}

\nc{\GCM}{GCM}
\nc{\HCO}{\widehat{\CO}}
\nc{\tCO}{\widetilde{\CO}}
\nc{\red}[1]{{\color{red}{#1}}}
\nc{\T}{\mathbb{T}}
\nc{\HC}{\mathsf{HC}}
\nc{\EV}{\mathsf{ev}}
\nc{\HRC}[1][n]{\widehat{\mathrm{RC}}_{{#1}}}
\nc{\hc}[1][n]{\ol{\mathcal{HC}}_{{#1}}}
\nc{\bphi}{\bar{\phi}}
\nc{\hgt}{\mathrm{ht}}

\begin{document}
\title{Graded Cartan determinants of the symmetric groups}
\author{Shunsuke Tsuchioka}
\address{
Kavli Institute for the Physics and Mathematics of the Universe, 
Todai Institutes for Advanced Study, 
the University of Tokyo, Kashiwa, Japan 277-8583 (Kavli IPMU)}
\thanks{The research was supported by
Grant-in-Aid for Research Activity Startup 22840026 and
Research Fellowships for Young Scientists 23$\cdot$8363,
Japan Society for the Promotion of Science.}
\email{tshun@kurims.kyoto-u.ac.jp}

\date{May 1, 2012}

\begin{abstract}
We give the graded Cartan determinants of the symmetric groups. 
Based on it, we propose a gradation of Hill's conjecture which is equivalent to K\"ulshammer-Olsson-Robinson's conjecture
on the generalized Cartan invariants of the symmetric groups.
\end{abstract}

\maketitle
%\tableofcontents
\section{Introduction}
Let $\Ggz$ be a finite-dimensional simple simply-laced complex Lie algebra.
The purpose of this paper is to give the Shapovalov determinants of the basic representation $V(\Lambda_0)$ of 
the untwisted quantum affinization $U_v(\Ggg)$ weight space wise (see 
Theorem \ref{MainTheorem}) and apply them to the graded representation theory of the symmetric groups (see \S\ref{GRC} and \S\ref{GRI}).

\subsection{Kashiwara's problem on the specialization for quantum groups}\label{KaPro}
The motivation of this paper comes from a problem of Kashiwara.
Let $\Gg$ be a symmetrizable Kac-Moody Lie algebra and let $\lambda\in \MP^+$ be a dominant integral weight.
In ~\cite[Problem 2]{Kas}, Kashiwara asks at which specialization $v=\xi\in\C^{\times}$ (see Remark \ref{speci})
the integrable highest weight $U_v(\Gg)$-module $V(\lambda)$ remains irreducible.
%Lusztig proved that the specialized module always remains irreducible
%if $\xi$ is not a root of unity.
When $\Gg$ is a finite-dimensional simple complex Lie algebra,
it is known that the specialized module remains irreducible
if $\xi$ is not a root of unity~\cite[Theorem 6.4]{APW}. More 
precisely, for an $\ell$-th root of unity $\xi=\exp(\frac{2\pi\sqrt{-1}}{\ell})$ where $\ell\geq 1$,
the specialized $U_v(\Gg)|_{v=\xi}$-module $V(\lambda)|_{v=\xi}$ is irreducible if $\ell<(\alpha,\lambda)$ for any 
positive root $\alpha\in\Delta^+$ of $\Gg$~\cite[Theorem 6.7]{APW}. 

When $\Gg$ is an affine Kac-Moody Lie algebra, almost nothing is known about Kashiwara's problem.
Chari-Jing proved that the specialized $U_v(\Ggg)|_{v=\xi}$-module 
$V(\Lambda_0)|_{v=\xi}$ is irreducible if $\xi$ is not a root of unity.
More precisely, for an $\ell$-th root of unity $\xi=\exp(\frac{2\pi\sqrt{-1}}{\ell})$ where $\ell\geq 1$,
if $\ell$ is coprime to
the Coxeter number $N$ of $\Ggz$, then the specialized module is irreducible~\cite[Theorem 4]{CJ}.
%$U_v(\Ggg)|_{v=\xi}$-module 
%$V(\Lambda_0)|_{v=\xi}$ is irreducible~\cite[Theorem 4]{CJ}.
Recall that $N=n$ (resp. $2(m-1),12,18,30$) when $\Ggz$ is of type $A_{n-1}$ for 
$n\geq 2$ (resp. $D_m$ for $m\geq 4$, $E_6$, $E_7$, $E_8$).
%As a corollary of the main result of this paper (Theorem \ref{MainTheorem}), we see that
The main result of this paper (Theorem \ref{MainTheorem}) gives an answer to Kashiwara's problem
when $\Gg=\Ggg$ and $\lambda$ is level 1.

\Th\label{kashiwara1}
Let $\Ggz$ be a finite-dimensional simple simply-laced complex Lie algebra of type $X$.
Then, for $\ell\geq 1$ the specialized $U_v(\Ggg)|_{v=\exp(\frac{2\pi\sqrt{-1}}{\ell})}$-module $V(\Lambda_0)|_{v=\exp(\frac{2\pi\sqrt{-1}}{\ell})}$ is 
irreducible if and only if 
\bna
\item $\GCD(\ell,2n)=1,2$ when $X=A_{n-1}$ for $n\geq 2$,
%\item $\GCD(\ell,n')=1$ and $\ell\not\in 4\Z$ when $X=A_{n-1}$ for even $n=2^{e}n'\geq 2$ where $e\geq 1$ and $n'\not\in 2\Z$,
\item $\ell\not\in 4\Z$ (resp. $\not\in 3\Z,\not\in 4\Z,\not\in 60\Z$) when $X=D_{n}$ for $n\geq 4$ (resp. $E_6,E_7,E_8$).
\ee
\enth

\subsection{Graded Cartan determinants}
Recently, Khovanov-Lauda and Rouquier independently
introduced a remarkable family of algebras (the KLR algebras) that
categorifies the negative half of the quantized enveloping algebras
associated with
symmetrizable Kac-Moody Lie algebras~\cite{KL1,KL2,Rou} (see Definition \ref{defKLR}).
An application of the KLR algebras is the
homogeneous presentation 
of the symmetric group algebras~\cite{BK2,Rou} (see Theorem \ref{BKR})
%of the symmetric group algebras and the Iwahori Hecke algebras of type A~\cite{BK1,Rou} (see Theorem \ref{BKR})
which quantize~\cite{BK1} Ariki's categorification of the Kostant $\Z$-form of the basic $\widehat{\mathfrak{sl}}_p$-module $V(\Lambda_0)^{\Z}\cong
\bigoplus_{n\geq 0}\KKK(\PROJ{\F_p\Sym_{n}})$.
Since the Shapovalov-like form (see Proposition \ref{exi_sh}) on the left hand side and the graded
Cartan pairing on the right hand side are compatible~\cite{BK3} (see Theorem ~\ref{compatibilityBK}),
the main result of this paper (Theorem \ref{MainTheorem}) can be interpreted as
the graded Cartan determinants of the symmetric groups.
The story is also valid for its $q$-analog, the Iwahori-Hecke algebra of type A (see \S\ref{IwahoriHecke}).

\Th\label{gradeddeter}
Let $p\geq 2$ be a prime number (resp. $p\geq 2$ an arbitrary integer).
Then, the graded Cartan determinant of a block of the symmetric group algebra over characteristic $p$ (resp. the Iwahori-Hecke algebra of type A over quantum characteristic $p$) 
%of $\mathbb{F}_p\mathfrak{S}_n$ (resp. $\mathcal{H}_n(\C;\exp(\frac{2\pi\sqrt{-1}}{p}))$) 
whose $p$-weight is $d$
is given by $\prod_{s=1}^{d}[p]_s^{N_{p,d,s}}$ where 
%(see also Notations in the end of this section)
\begin{align*}
N_{p,d,s}=\sum_{\lambda\in\PAR(d)}\frac{m_s(\lambda)}{p-1}\prod_{u\geq 1}\binom{m_u(\lambda)+p-2}{m_u(\lambda)}
=\sum_{(\lambda_i)_{i=1}^{p-1}\in\PAR_{p-1}(d)}m_s(\lambda_1).
\end{align*}
For the notations on partitions, see Notations in the end of \S1.
\enth

In virtue of the classical results of the modular representation theory of finite groups,
the Cartan determinant (see Definition \ref{CartanInv}) of (a block of) a group algebra of a finite group $G$ over characteristic $p\geq 2$
is a power of $p$~\cite[Theorem 1]{Bra} and its degree reflects $p$-defects of $p$-regular conjugacy classes of $G$~\cite[Part III,\S16]{BrNe} (see 
Theorem \ref{classicalmodular}).
Based on Olsson's calculation of $p$-defects of the symmetric groups~\cite{Ols},
Bessenrodt-Olsson gave an explicit formula for the Cartan determinants of (a block of) a symmetric group algebra
over characteristic $p$ in terms of generating functions~\cite[Theorem 3.3, Theorem 3.4]{BO1}.

On the other hand, Brundan-Kleshchev arrived at the same formula %as Bessenrodt-Olsson 
by a different method~\cite[Corollary 1]{BK4}.
They utilized Ariki's categorification and
turned the calculation of the $p$-Cartan determinants of the symmetric groups %over characteristic $p$ 
into the calculation of the Shapovalov determinants of the basic $\widehat{\mathfrak{sl}}_p$-module $V(\Lambda_0)$.
Their block-theory-free method has an advantage in that
it is applicable to the Iwahori-Hecke algebra of type A.
In fact, Brundan-Kleshchev settled affirmatively the conjecture of Mathas
that the Cartan determinants of (a block of) the Iwahori-Hecke algebra of type A over
quantum characteristic $\ell$ is a power of $\ell$ for arbitrary $\ell\geq 2$ (see ~\cite{Don}).

Theorem \ref{gradeddeter} (and Theorem \ref{MainTheorem}) can be seen a gradation of Brundan-Kleshchev's result
and our method is similar to theirs. 
A difference is that we use a bi-additive form that obeys
a rule slightly different from that of the Shapovalov form (see Proposition \ref{exi_sh}).

%As we have said, both Theorem \ref{kashiwara1} and Theorem \ref{gradeddeter} follow from Theorem \ref{MainTheorem}.
%Theorem \ref{MainTheorem} gives the Shapovalov determinants of the basic representations $V(\Lambda_0)$ of the quantum groups $U_v(\Ggg)$.
%Brundan-Kleshchev calculated the Shapovalov determinants at $v=1$~\cite[Main Theorem]{BK4}.
%Our method is similar to theirs. 
%A difference is that we use a bi-additive form (which we call quasi-Shapovalov form) that obeys
%a rule slightly different from that of the Shapovalov form (see Proposition \ref{exi_sh}).

We remark that Brundan-Kleshchev also calculated the Shapovalov determinants of the basic %level 1 integrable highest weight 
module
over a twisted affinization of 
%a finite-dimensional semisimple simply-laced Lie algebra 
$\Ggz$. % at $v=1$.
Although our method should be applicable in the quantum setting,
%calculation of the Shapovalov determinants of the basic
%module
%over a quantum twisted affinization of $\Ggz$ 
%in principle,
there are some obstructions to overcome (see Remark \ref{Drinfeldtwisted}). 
We discuss them and give a conjectural Shapovalov determinants (see Conjecture \ref{conjshap}).

Recently, the authors~\cite{KKT} introduced a superversion of the KLR algebras which they call quiver Hecke superalgebras
and showed ``Morita superequivalence'' between the spin symmetric group algebras %(or its quantization the Hecke-Clifford superalgebras)
and quiver Hecke superalgebras~\cite[Theorem 5.4]{KKT}.
From the results of ~\cite{BK5,Tsu}, 
it is expected that the quiver Hecke superalgebras categorify the quantum groups (see also ~\cite{HW}) and through this
the Shapovalov determinants of the basic module of $U_v(A^{(2)}_{p-1})$ for odd prime $p\geq 3$
coincides with the graded Cartan determinants of the spin symmetric group algebras. %(see Conjecture \ref{conjspincartan}).

\subsection{K\"ulshammer-Olsson-Robinson's conjecture and graded Cartan invariants}\label{introinv}
Based on an observation that combinatorial notions for partitions (such as $p$-cores) that appear in the description 
of the $p$-modular representation-theoretic invariants of the symmetric groups
make sense when $p$ is not a prime, K\"ulshammer-Olsson-Robinson initiated the study of 
``$\ell$-modular representations of the symmetric groups'' for arbitrary $\ell\geq 2$~\cite{KOR}.

Their main result is a generalization of Nakayama-Brauer-Robinson classification of
the $p$-blocks of $\mathfrak{S}_n$ for an any $\ell\geq 2$~\cite[Theorem 5.13]{KOR} (see Theorem ~\ref{KORmainthm}).
Based on it, they gave a conjectural generalization of the Cartan invariants of the symmetric groups
for an arbitrary $\ell\geq 2$~\cite[Conjecture 6.4]{KOR} (see Conjecture ~\ref{KORcon} which we call the KOR conjecture for short). 
As an evidence, they checked that their invariants gives a generalized Cartan determinant
calculated by Brundan-Kleshchev~\cite{BK4}.

Hill settled affirmatively the KOR conjecture for any $\ell\geq 2$ such that each prime factor $p$ of $\ell$ divides $\ell$ at most $p$ times~\cite[Theorem 1.3]{Hil}.
In his course of proof, Hill gave a conjecture~\cite[Conjecture 10.5]{Hil} (see Conjecture \ref{Hcon}).
Hill's conjecture and the KOR conjecture are equivalent in virtue of ~\cite{BH} (see Corollary \ref{equiva}).

As graded Cartan determinants %(Theorem \ref{gradeddeter2}) 
calculated in this paper,
it is reasonable to expect that a proof of Hill's conjecture similarly works in a graded setting.
% and
%that the correct gradation of Hill's conjecture reveals some essence of the original conjecture
%which is hidden in the ungraded setting. 
Motivated by this, we propose a gradation of Hill's conjecture (see Conjecture \ref{myconj}).
%Conjecture \ref{myconj} is drawn from computer calculation and holds at least for $pr\leq 9$ and $d\leq 15$ (see Remark \ref{computer_check}).
%When $r=1$, Conjecture \ref{myconj} is the same as ~\cite[Conjecture 8.2 (i)]{ASY}.
As a support, we check that Conjecture \ref{myconj} is compatible with the graded Cartan determinants (see Theorem \ref{conjcheck}).
%of Theorem \ref{gradeddeter2} (see Theorem \ref{conjcheck}).
The proof itself may be of interest.
%We expect that our gradation is correct and conjecture reveals some essence of the original conjecture
%which is hidden in the ungraded setting. 
We expect that our gradation is correct and gives an insight to future trials of the proof
of Hill's conjecture (and thus the KOR conjecture).
%of Hill's conjecture and thus contribute to the proof of the KOR conjecture.

\vskip 3mm

\noindent{\bf Notations.} $\N=\Z_{\geq 0}$ (resp. $\NNN=\Z_{\geq 1}$) means the set of non-negative (resp. positive) integers.
We denote the empty partition by $\emptypartition$ and reserve the symbol $\emptyset$ for the empty set.
For a partition $\lambda=(\lambda_1,\lambda_2,\cdots)$, we define $m_k(\lambda)=|\{i\geq 1\mid \lambda_i=k\}|$ for $k\geq 1$.
We also define $|\lambda|=\sum_{i\geq 1}\lambda_i$ and $\ell(\lambda)=\sum_{i\geq 1}m_i(\lambda)$.
We denote by $\PAR(n)$ the set of partitions of $n\geq 0$ and define $\PAR=\bigsqcup_{n\geq 0}\PAR(n)$.
For $m,n\geq 0$, we denote by $\PAR_m(n)$ the multipartitions of $n$, i.e., $\PAR_m(n)=\{(\lambda_i)_{i=1}^{m}\in\PAR^m\mid \sum_{i=1}^{m}|\lambda_i|=n\}$
and put $u(m,n)=|\PAR_m(n)|$. Note that $u(0,0)=1$ and $u(0,n)=0$ for $n\geq 1$. 
%For $m\geq 0$, the generating function is given by $\sum_{n\geq 0}u(m,n)x^n=1/\PHI(x)^m$ where $\PHI(x)=\prod_{k\geq 1}(1-x^k)$.

\vskip 3mm

\noindent{\bf Organization of the paper}
%The paper is organized as follows.
In \S2, we develop a linear algebra for Gram determinants.
In \S3, we review the boson-fermion correspondence %over $\Z$ due to DeConcini-Kac-Kazhdan 
and deduce the existence of a bi-additive form we need.
In \S4, we give a brief review on quantum groups.
In \S5, we calculate the Shapovalov determinants of the basic representations of the untwisted affine A,D,E quantum groups
and in \S6 using the KLR algebras we interpret it as graded Cartan determinants of the symmetric groups and its Iwahori-Hecke akgebras.
Finally, in \S7 we propose a gradation on Hill's conjecture and make some justifications.

\vskip 3mm

\noindent{\bf Acknowledgements.} The content of \S7 is benefited from the discussion with Yasuhide Numata.
The author would like to thank him and also
%The author would also 
would like to thank Hiraku Nakajima and Yoshiyuki Kimura
for discussions on quantum groups. %The content of \S7 is benefited by the discussion with Yasuhide Numata. 
The author is grateful to Hiro-Fumi Yamada. 
Without his encouragement and their paper ~\cite{ASY}, this paper would not have been written.

\section{Linear Algebra on Gram determinants}
Let $\cor$ be a field and let $V$ be a $\cor$-vector space.
We say that a map $B:V\times V\to\cor$ is bi-additive if we have $B(w_1+w_2,w)=B(w_1,w)+B(w_2,w)$ and $B(w,w_1+w_2)=B(w,w_1)+B(w,w_2)$ for all
$w,w_1,w_2\in V$.
The radical $\RAD(B)$ of a bi-additive map $B$ is defined to be one of the following two additive subgroups of $V$ 
\begin{align*}
\{
w_1\in V\mid\forall w_2\in V, B(w_1,w_2)=0
\},\quad
\{
w_2\in V\mid\forall w_1\in V, B(w_1,w_2)=0
\}
\end{align*}
if they are equal (otherwise, we do not define $\RAD(B)$). 
We say that a bi-additive map $B$ is non-degenerate if $\RAD(B)$ is defined and $\RAD(B)=0$.

%\subsection{Gram matrices and Gram determinants}
\Def\label{GRM} 
%In this subsection,
Let $(\cor,\BAR,\MA,V,\BL,V^\MA)$ be a 6-tuple such that
\bna
\item $\cor$ is a field with a ring involution $\BAR:\cor\isoto\cor$,
\item $\MA$ is a subring of $\cor$ such that $\BAR(\MA)\subseteq \MA$,
\item $V$ is a finite dimensional $\cor$-vector space,
\item $V^\MA$ is an $\MA$-lattice of $V$, i.e., $V^\MA$ is a free $\MA$-module with an isomorphism $\cor\otimes_\MA V^\MA\isoto V$,
\item $\BL:V\times V\to \cor$ is a bi-additive non-degenerate map such that
$\BL(aw_1,w_2)=\BAR(a)\BL(w_1,w_2)$ and $\BL(w_1,aw_2)=a\BL(w_1,w_2)$ for all $w_1,w_2\in V$ and $a\in \cor$.
\ee

%Let $(w_i)_{1\leq i\leq \dim V}$ be an $\MA$-basis of $V^\MA$. 
To such a 6-tuple, we define
the Gram matrix $\GRM_{\cor,\BAR,\MA}(V,\BL,V^\MA)=(\BL(w_i,w_j))_{1\leq i,j\leq \dim V}$ and
define the Gram determinant $\GR_{\cor,\BAR,\MA}(V,\BL,V^\MA)=\det\GRM_{\cor,\BAR,\MA}(V,\BL,V^\MA)$
where $(w_i)_{1\leq i\leq \dim V}$ is an $\MA$-basis of $V^\MA$.
\edf

Let $X$ and $Y$ are two Gram matrices associated with two choices of $\MA$-bases of $V^\MA$. % of the 6-tuple of Definition \ref{GRM}.
Clearly, there exists $P\in \GL_{\dim V}(\MA)$ such that $X=\BAR(\TRANS{P})YP$.
Thus, $\GR_{\cor,\BAR,\MA}(V,\BL,V^\MA)$ is uniquely determined as an element of $\cor^\times/\{a\BAR(a)\mid a\in \MA^\times\}$.

%\Def
%Let $(\cor,\BAR,\MA,V,\BL,V^\MA)$ be a 6-tuple as in Definition \ref{GRM}.
%We define the Gram determinant $\GR_{\cor,\BAR,\MA}(V,\BL,V^\MA)$ as 
%$\det\GRM_{\cor,\BAR,\MA}(V,\BL,V^\MA)\in\cor^\times/\{a\BAR(a)\mid a\in \MA^\times\}$.
%\edf

%\subsection{A comparison result}
The following comparison result is based on ~\cite[\S5]{BK4}.

\Prop\label{comp} Let $I$ be a finite set and let $\cor$ be a field. % of characteristic 0 with a ring involution $\BAR:\cor\isoto\cor$.
We regard a polynomial ring $V=\cor[y_n^{(i)}\mid i\in I,n\geq 1]$ as a graded $\cor$-algebra via $\DEG y_n^{(i)}=n$ and
denote by $V_d$ the $\cor$-vector subspace of $V$ consisting of homogeneous elements of degree $d$ for $d\geq 0$.
%Further, 
Assume $\CHAR\cor=0$ and we are given the following data.%\footnote{This assumption automatically implies that $\CHAR\cor=0$.}.
\bna
\item a subring $\MA$ of $\cor$ and a ring involution $\BAR:\cor\isoto\cor$ such that $\BAR(\MA)\subseteq \MA$,
%\item a finite set $I$,
\item a family of invertible matrices $A=(A^{(m)})_{m\geq 1}$ where $A^{(m)}=(a^{(m)}_{ij})_{i,j\in I}\in\GL_I(\MA^\BAR)$,
\item(\label{lin_change2}) two bi-additive non-degenerate maps $\langle,\rangle_S, \langle,\rangle_K:V\times V\to \cor$ 
such that
\begin{itemize}
\item $\langle af,g\rangle_X=\BAR(a)\langle f,g\rangle_X$, $\langle f,ag\rangle_X=a\langle f,g\rangle_X$ and 
$\langle f,g\rangle_X=\BAR(\langle g,f\rangle_X)$
\item $\langle 1,1\rangle_S=\langle 1,1\rangle_K$ and 
$\langle my_m^{(i)}f,g\rangle_S=\langle f,\sum_{j\in I}a_{ij}^{(m)}\frac{\partial g}{\partial y^{(j)}_m}\rangle_S$, 
$\langle my_m^{(i)}f,g\rangle_K=\langle f,\frac{\partial g}{\partial y^{(i)}_m}\rangle_K$,% for all $f,g\in V$ and $m\geq 0$,
\end{itemize}
for $X\in\{S,K\}$ and $f,g\in V,a\in\cor,m\geq1$, % (such maps exist only when $\CHAR\cor=0$),
\item\label{lin_change} a family of new variables $(x_n^{(i)})_{\substack{i\in I, \\ n\geq 1}}$ such 
that $x_n^{(i)}-y_n^{(i)}\in \cor[y_m^{(i)}\mid 1\leq m<n]\cap V_n$.% for $i\in I,n\geq 1$.
\ee
Then, in $\cor^\times/\{a\BAR(a)\mid a\in \MA^\times\}$ we have
%for any $d\geq 0$
\begin{align*}
\GR_{\cor,\BAR,\MA}(V_d,\langle,\rangle_S,V_d^\MA)/\GR_{\cor,\BAR,\MA}(V_d,\langle,\rangle_K,V_d^\MA)=\Delta_d(A)
\end{align*}
for any $d\geq 0$
where $V_d^\MA = V_d\cap \MA[x_n^{(i)}\mid i\in I,n\geq 1]$ and
\begin{align*}
\Delta_d(A) = \prod_{s=1}^{d}(\det A^{(s)})^{\sum_{\lambda\in\PAR(d)}\frac{m_s(\lambda)}{|I|}\prod_{u\geq 1}\binom{m_u(\lambda)+|I|-1}{m_u(\lambda)}}.
\end{align*}
\enprop

\Proof Firstly, we note that in virtue of (\ref{lin_change}) we have $V=\cor[x_n^{(i)}\mid i\in I,n\geq 1]$.
Thus, $V_d^\MA$ is an $\MA$-lattice of $V_d$.

From now on, we fix a total ordering on $I$ and define
\begin{align*}
\Omega(\lambda) = \{
(i_1,\cdots,i_{\ell(\lambda)})\in I^{\ell(\lambda)}\mid \lambda_j=\lambda_{j+1}\Rightarrow i_j\geq i_{j+1}
\}%,\quad
%\Omega_d = \{(\lambda,i)\mid \lambda\in\PAR(d),i\in \Omega(\lambda)\}
\end{align*}
for $\lambda\in\PAR$ and put $\Omega_d = \{(\lambda,i)\mid \lambda\in\PAR(d),i\in \Omega(\lambda)\}$ for $d\geq 0$.

We define a family of new variables $z^{(i)}_n=\sum_{j\in I}a_{ij}^{(n)}y^{(j)}_n$ for $i\in I$ and $n\geq 1$.
Then, we see that for any $d\geq 0$
\begin{itemize}
\item $\{x^{(i)}_{\lambda}:=\prod_{k=1}^{\ell(\lambda)}x^{(i_k)}_{\lambda_k}\}_{(\lambda,i)\in\Omega_d}$
is an $\MA$-basis of $V^\MA_d$,
\item $\{y^{(i)}_{\lambda}:=\prod_{k=1}^{\ell(\lambda)}y^{(i_k)}_{\lambda_k}\}_{(\lambda,i)\in\Omega_d}$ and 
$\{z^{(i)}_{\lambda}:=\prod_{k=1}^{\ell(\lambda)}z^{(i_k)}_{\lambda_k}\}_{(\lambda,i)\in\Omega_d}$
are $\cor$-basis of $V_d$.
\end{itemize}

We show for any $\lambda\in\PAR$ by induction on $\ell:=\ell(\lambda)\geq 0$ that we have
\begin{align}
\forall i\in\Omega(\lambda),
\forall f\in V,\langle y^{(i)}_{\lambda}, f \rangle_S=\langle z^{(i)}_{\lambda}, f \rangle_K.
\label{henkan_yz}
\end{align}
%for any $\lambda=(\lambda_1,\cdots,\lambda_{\ell})\vdash d, i\in\Omega(\lambda)$ 
It is clear by (\ref{lin_change2}) when $\ell=0$ and 
assume that (\ref{henkan_yz}) holds for any $\lambda$ with $\ell(\lambda)<\ell$.
Put $j=(i_2,\cdots,i_\ell)$ and $\mu=(\lambda_2,\cdots,\lambda_\ell)$ and (\ref{henkan_yz}) is proved as follows.
\begin{align*}
\langle y^{(i)}_{\lambda}, f \rangle_S =\langle y^{(i_1)}_{\lambda_1}y^{(j)}_{\mu}, f \rangle_S 
&= \frac{1}{\lambda_1}\langle y^{(j)}_{\mu}, \sum_{k\in I}a_{i_1,k}^{(\lambda_1)}\frac{\partial f}{\partial y^{(k)}_{\lambda_1}} \rangle_S \\
& =\frac{1}{\lambda_1}\langle z^{(j)}_{\mu}, \sum_{k\in I}a_{i_1,k}^{(\lambda_1)}\frac{\partial f}{\partial y^{(k)}_{\lambda_1}} \rangle_K 
= \langle z^{(i_1)}_{\lambda_1}z^{(j)}_{\mu}, f \rangle_K =\langle z^{(i)}_{\lambda}, f \rangle_K.
\end{align*}

We note that $\GR_{\cor,\BAR,\MA}(V_d,\langle,\rangle_S,V_d^\MA)=\det M_S$ and $\GR_{\cor,\BAR,\MA}(V_d,\langle,\rangle_K,V_d^\MA)=\det M_K$ where
$M_S=(\langle x^{(i)}_{\lambda},x^{(j)}_{\mu}\rangle_S)_{(\lambda,i),(\mu,j)\in\Omega_d}$ and
$M_K=(\langle x^{(i)}_{\lambda},x^{(j)}_{\mu}\rangle_K)_{(\lambda,i),(\mu,j)\in\Omega_d}$.

Let $P\in \GL_{\Omega_d}(\cor)$ and $Q\in \GL_{\Omega_d}(\MA^\BAR)$ by
\begin{align*}
x_{\lambda}^{(i)}=\sum_{(\mu,j)\in\Omega_d}p^{(i,j)}_{\lambda,\mu}y^{(j)}_{\mu},\quad
z_{\lambda}^{(i)}=\sum_{(\mu,j)\in\Omega_d}q^{(i,j)}_{\lambda,\mu}y^{(j)}_{\mu}.
\end{align*}
In virtue of (\ref{henkan_yz}), we have 
\begin{align*}
M_S &= 
\BAR(P)(\langle y^{(i)}_{\lambda},x^{(j)}_{\mu}\rangle_S)_{(\lambda,i),(\mu,j)\in\Omega_d}=
\BAR(P)(\langle z^{(i)}_{\lambda},x^{(j)}_{\mu}\rangle_K)_{(\lambda,i),(\mu,j)\in\Omega_d}, \\
M_K &= 
\BAR(P)(\langle y^{(i)}_{\lambda},x^{(j)}_{\mu}\rangle_K)_{(\lambda,i),(\mu,j)\in\Omega_d}=
\BAR(P)Q^{-1}(\langle z^{(i)}_{\lambda},x^{(j)}_{\mu}\rangle_K)_{(\lambda,i),(\mu,j)\in\Omega_d}.
\end{align*}
Thus, we have $M_S=\BAR(P)Q\BAR(P)^{-1}M_K$ and it is enough to show that $\det Q=\Delta_d$.

It is easily seen that $q^{(i,j)}_{\lambda,\mu}=0$ for $(\lambda,i),(\mu,j)\in\Omega_d$ with $\lambda\ne\mu$.
Since $(q^{(i,j)}_{\lambda,\lambda})_{i,j\in\Omega(\lambda)}$ is identified with $\otimes_{t\geq1}\SYM^{m_t(\lambda)}(A^{(t)})$,
we have
\begin{align*}
\det Q = \prod_{\lambda\in\PAR(d)}\det ((q^{(i,j)}_{\lambda,\lambda})_{i,j\in\Omega(\lambda)}) 
= \prod_{\lambda\in\PAR(d)}\det(\otimes_{t\geq1}\SYM^{m_t(\lambda)}(A^{(t)})).
\end{align*}

Apply Lemma \ref{easy_lin_alg} which is an easy linear algebra exercise, we get
\begin{align*}
\det Q &= \prod_{\lambda\in\PAR(d)}\prod_{t\geq1}\det(\SYM^{m_t(\lambda)}(A^{(t)}))^{\prod_{u\ne t}\dim \SYM^{m_u(\lambda)}(k^{I})} \\
&= \prod_{\lambda\in\PAR(d)}\prod_{t\geq1}(\det(A^{(t)}))^{\binom{|I|+m_t(\lambda)-1}{m_t(\lambda)-1}\prod_{u\ne t}\binom{|I|+m_u(\lambda)-1}{m_u(\lambda)}}.
\end{align*}

Since $\binom{|I|+m_t(\lambda)-1}{m_t(\lambda)-1}\prod_{u\ne t}\binom{|I|+m_u(\lambda)-1}{m_u(\lambda)}=\frac{m_t(\lambda)}{|I|}\prod_{u}\binom{|I|+m_u(\lambda)-1}{m_u(\lambda)}$, we have $\det Q=\Delta_d$.
%&= \prod_{\lambda\in\PAR(d)}\prod_{t=1}^{s}(\det(A^{(t)}))^{\frac{m_t(\lambda)}{|I|}\prod_{u=1}^{s}\binom{|I|+m_u(\lambda)-1}{m_u(\lambda)}} = \Delta_d(A).
\QED

\Lemma\label{easy_lin_alg}
Let $\cor$ be a field.
\bna
\item Let $V$ be an $n$-dimensional 
$\cor$-vector space and let $f\in\END_\cor(V)$. For $m\geq 0$, we have
\begin{align*}
\dim\SYM^mV=\binom{n+m-1}{m},\quad
\det\SYM^mf=(\det f)^{\binom{n+m-1}{m-1}}.
\end{align*}
\item For an $\ell$-tuple of matrices $(B^{(i)})_{1\leq i\leq \ell}$ with $B^{(i)}\in\MAT_{n_i}(\cor)$, 
we have $\det(\otimes_{i=1}^{\ell}B^{(i)})=\prod_{i=1}^{\ell}(\det B^{(i)})^{N/n_i}$ %for $\{B^{(i)}\in\MAT_{n_i}(\cor)\mid 1\leq i\leq \ell\}$
where $N=\prod_{i=1}^{\ell}n_i$.
\ee
\enlemma

\section{Boson-fermion correspondence}

%Throughout this section, let $R$ be a commutative ring with $1$.

\Def Let $\MAA$ be a commutative ring with $1$.
\bna
\item The fermionic Fock space $\FF_\MAA=\MAA^{\oplus\SMT}$ over $\MAA$ is a free $\MAA$-module with a basis consisting of all semi-infinite monomials
\begin{align*}
\SMT:=\{i_1 \wedge i_2 \wedge\cdots\mid \textrm{$\forall k\geq 1,i_k>i_{k+1}$ and $\exists N>0,\forall k>N,i_k=i_{k-1}-1$}\}.
\end{align*}
\item The bosonic Fock space $\BF_\MAA$ over $\MAA$ is a polynomial ring $\MAA[v,v^{-1},y_1,y_2,\cdots]$.
\ee
\edf

For $m\in\Z$, we denote by $|m\rangle$ the vacuum vector of charge $m$, i.e., $|m\rangle=m \wedge (m-1)\wedge\cdots\in\SMT$.
We say that a semi-infinite monomial $i_1 \wedge i_2\wedge\cdots\in\SMT$ is of charge $m$ if
it differs from $|m\rangle$ only at a finite number of places, i.e.,
there exists $N>0$ such that we have $i_k=m+1-k$ for all $k>N$.

Recall that for $j\in\Z$
a creation operator $\psi_j\in\END_\MAA(\FF_\MAA)$ and an annihilation operator $\psi^\ast_j\in\END_\MAA(\FF_\MAA)$ are defined as follows
where we regard $i_0$ as $+\infty$.
\begin{align*}
\psi_j(i_1 \wedge i_2 \wedge\cdots)
&=
\begin{cases}
0 & \exists s\geq1,j=i_s, \\
(-1)^s i_1\wedge \cdots \wedge i_s \wedge j \wedge i_{s+1}\wedge\cdots & \exists s\geq0,i_s>j>i_{s+1}, 
\end{cases} \\
\psi^\ast_j(i_1 \wedge i_2 \wedge\cdots)
&=
\begin{cases}
0 & \forall s\geq1,j\ne i_s, \\
(-1)^{s+1} i_1\wedge \cdots \wedge i_{s-1} \wedge i_{s+1}\wedge\cdots & \exists s\geq1,i_s=j.
\end{cases}
\end{align*}

\Rem\label{biadjoint}
For any $s,t\in\SMT$ and $j\in\Z$, we have $\psi_j(s)=t$ if and only if $\psi^{\ast}_j(t)=s$.
\enrem

\Rem\label{nization}
$\BF_\MAA$ has a natural structure of module over the 
Heisenberg algebra $\MAAA_\MAA:=\langle\{\alpha_m\}_{m\in\Z}\mid \alpha_m\alpha_n-\alpha_n\alpha_m=m\delta_{m+n,0}\rangle_{\textrm{$\MAA$-alg}}$
through the assignment of linear operators %where $k\geq 1$
\begin{align*}
\alpha_0=v\frac{\partial}{\partial v},\quad
\alpha_k=\frac{\partial}{\partial y_k},\quad
\alpha_{-k}=ky_k
\end{align*}
where $k\geq 1$.
On the other hand, bosonization turns $\FF_\MAA$ into an $\MAAA_\MAA$-module where $\ell\ne 0$
\begin{align*}
\alpha_k = \sum_{j\in\Z}\psi_j\psi^\ast_{j+k},\quad
\alpha_0 = \sum_{j>0}\psi_j\psi^\ast_j-\sum_{j\leq 0}\psi^\ast_j\psi_j.
\end{align*}
\enrem

\Def
For a partition $\lambda\in\PAR$,
the Schur function $s_{\lambda}\in\Q[y_1,y_2,\cdots]$ is defined as 
$s_{\lambda}=\det(x_{\lambda_i-i+j})_{1\leq i,j\leq|\lambda|}$ where $x_{m}=\delta_{m,0}$ for $m\leq 0$ and
\begin{align}
1+\sum_{n\geq 1}x_nz^n = \exp(\sum_{n\geq 1}y_nz^n).
\label{xygenseries}
\end{align}
\edf

\Prop[{\cite[\S14.10]{Kac}}] Let $\cor$ be a field of characteristic 0.
\bnum
\item There exists a unique $\MAAA_\cor$-module isomorphism $\sigma_\cor:\FF_\cor\isoto\BF_\cor$ such that
$\sigma_\cor(|m\rangle)=v^m$ for any $m\in\Z$.
\item For any $m\in\Z$ and $\WI=i_1\wedge i_2\wedge\cdots\in\SMT$ of charge $m$, 
we have $\sigma_\cor(\WI)=v^ms_{\lambda(\WI)}$ where $\lambda(\WI)=(i_1-m,i_2-(m-1),\cdots)\in\PAR$.
\ee
\enprop

\Cor[{\cite[Theorem 2]{DcKK}}]
Let $\cor$ be a field of characteristic 0.
Define a family of new variables $(x_n)_{n\geq 1}$ by (\ref{xygenseries})
and put
$\BFF_\Z=\Z[v,v^{-1},x_1,x_2,\cdots](\subseteq \BF_\cor)$. %=\BFF_\cor:=\cor[v,v^{-1},x_1,x_2,\cdots]).
Then, the restriction $\sigma_\Z=\sigma_\cor|_{\FF_\Z}$ induces an $\MAAA_{\Z}$-module isomorphism 
$\sigma_\Z:\FF_\Z\isoto \BFF_\Z$.
\encor

\Cor
Let $\cor$ be a field. % of characteristic 0. % with a ring involution $\BAR:\cor\isoto\cor$. 
We regard a polynomial ring $V=\cor[y_1,y_2,\cdots]$ as a graded $\cor$-algebra via $\DEG y_n=n$ and
denote by $V_d$ the $\cor$-vector subspace of $V$ consisting of homogeneous elements of degree $d$ for $d\geq 0$.
Assume $\CHAR\cor=0$ and take a ring involution $\BAR:\cor\isoto\cor$.% be a ring involution.
\bna
\item There
exists a unique bi-additive non-degenerate map $\langle,\rangle_K:V\times V\to\cor$ such that
\bnum
\item $\langle af,g\rangle_K=\BAR(a)\langle f,g\rangle_K, \langle f,ag\rangle_K=a\langle f,g\rangle_K$ % for any $f,g\in V$ and $a\in \cor$,
and $\langle f,g\rangle_K=\BAR(\langle g,f\rangle_K)$, % for any $f,g\in V$,
\item $\langle 1,1\rangle_K=1$ and 
$\langle my_m^{(i)}f,g\rangle_K=\langle f,\frac{\partial g}{\partial y^{(i)}_m}\rangle_K$, %for any $f,g\in V$.
\ee
for any $f,g\in V$ and $a\in \cor$.
\item Let $\MA$ be a subring of $\cor$ such that $\BAR(\MA)\subseteq\MA$.
We choose an $\MA$-lattice $V^\MA_d$ of $V_d$ by $V^\MA_d=V_d\cap\MA[x_n\mid n\geq 1]$ for $d\geq 0$ where
$(x_n)_{n\geq 1}$ is defined in terms of $(y_n)_{n\geq 1}$ via (\ref{xygenseries}).
Then, for all $d\geq 0$ we have $\GR_{\cor,\BAR,\MA}(V_d,\langle,\rangle_K,V^{\MA}_d)=1$ in $\cor^\times/\{a\BAR(a)\mid a\in\MA^{\times}\}$.
\ee
\encor

\Proof (a)\ Uniqueness is clear and we prove the existence.
We introduce a bi-additive non-degenerate map $(,):\FF_\cor\times \FF_\cor\to\cor$ as 
$(\sum_{s\in\SMT}a_s\cdot s,\sum_{t\in\SMT}b_t\cdot t)=\sum_{u\in\SMT}\BAR(a_u)b_u$.
Let $\FF^{(0)}_\cor$ be the $\cor$-subspace of $\FF_\cor$ spanned by $\{s\in\SMT\mid\textrm{$s$ is of charge $0$}\}$.
Then, $\sigma_{\cor}|_{\FF^{(0)}_\cor}:\FF^{(0)}_\cor\isoto V(\subseteq \BF_\cor)$ carries $(,)|_{\FF^{(0)}_\cor\times\FF^{(0)}_\cor}$ on $\FF^{(0)}_\cor$ 
into the desired map $\langle,\rangle_K$ on $V$ because of Remark \ref{biadjoint} and Remark \ref{nization}.

\smallskip
\noindent (b) It is well-known that $\{s_{\lambda}\mid \lambda\in\PAR\}$ forms a $\Z$-basis of $V^{\Z}:=\Z[x_1,x_2,\cdots]\subseteq V$.
From the construction of $\langle,\rangle_K$ above,
it is clear that $\{s_{\lambda}\mid \lambda\in\PAR\}$ is orthonormal on $V^{\Z}$ (and thus $\MA\otimes V^{\Z}$).
\QED

\Cor\label{WITHOUTK}
In the setting of Proposition \ref{comp}, we may take $x^{(i)}_n$ for $i\in I$ and $n\geq 1$ by
\begin{align*}
1+\sum_{n\geq 1}x^{(i)}_nz^n = \exp(\sum_{n\geq 1}y^{(i)}_nz^n).
\end{align*}
Then, we have $\GR_{\cor,\BAR,\MA}(V_d,\langle,\rangle_S,V_d^\MA)=\Delta_d(A)$.
\encor

\section{Quantum groups and Shapovalov forms}

Let $v$ be an indeterminate. In the rest of the paper, we work in the field $\cor=\Q(v)$ and its subring $\MA=\Z[v,v^{-1}]$.
We consider two ring involution of $\cor$. One is the identity map $\ID_{\cor}$ and the other is
a $\Q$-algebra involution $\BAR:\cor\to\cor$ which maps $v$ to $v^{-1}$.

\subsection{Quantum groups}

Let $A=(a_{ij})_{i,j\in I}$ be a symmetrizable {\GCM} and take the symmetrization $d=(d_i)_{i\in I}$ of $A$, i.e.,
a unique $d\in\NNN^I$ such that $d_ia_{ij}=d_ja_{ji}$ for all $i,j\in I$ and
$\gcd(d_i)_{i\in I}=1$.
We take a Cartan data $(\MP,\MPC,\Pi,\Pi^\vee)$ in the following sense.
\bna
\item $\MPC$ is a free $\Z$-module of rank $(2|I|-\Rank A)$ and $\MP=\Hom_\Z(\MPC,\Z)$,
\item $\Pi^\vee=\{h_i\mid i\in I\}$ is a $\Z$-linearly independent subset of $\MPC$,
\item $\Pi=\{\alpha_i\mid i\in I\}$ is a $\Z$-linearly independent subset of $\MP$,
\item $\alpha_j(h_i)=a_{ij}$ for all $i,j\in I$.
\ee

We denote by $\MP^+$
the set of dominant integral weights $\{\lambda\in\MP\mid \forall i\in I,\lambda(h_i)\in\N\}$.
For each $i\in I$, $\Lambda_i\in\MP^+$ is a dominant integral weight 
determined modulo the subgroup $\{\lambda\in\MP\mid \forall i\in I,\lambda(h_i)=0\}(\subseteq \MP)$ 
by the condition $\Lambda_i(h_j)=\delta_{ij}$ for all $j\in I$.

Recall that
the Weyl group $W=W(A)$ is a subgroup of $\GL(\Gh^\ast)$ generated by $\{s_i:\Gh^\ast\isoto\Gh^\ast,\lambda\MAPSTO\lambda-\lambda(h_i)\alpha_i\mid i\in I\}$ and
that $W$ also acts on $\MP(\subseteq\Gh^\ast)$.
% where
%$
%s_i:\Gh^*\isoto\Gh^*,%\quad
%\lambda\MAPSTO\lambda-\lambda(h_i)\alpha_i.
%$

In the following, we use the usual abberiviation such as 
$v_i=v^{d_i}$, $[n]_{\ell}=\sum_{k=1}^{n}v^{(n+1-2k)\ell}$, $[n]_{\ell}!=\prod_{m=1}^{n}[m]_{\ell}$
and $\Binom{n}{m}{i}=\frac{[n]_i!}{[m]_i![n-m]_i!}$ 
for $i\in I$ and $n\geq m\geq 0$.

\Def 
The quantum group $U_v=U_v(A)$ is an associative $\cor$-algebra with $1$ generated by 
$\{e_i,f_i\mid i\in I\}\cup\{v^h\mid h\in\MPC\}$ with the following definition relations.
\bna
\item $v^0=1$ and $v^hv^{h'}=v^{h+h'}$ for any $h,h'\in\MP$,
\item $v^{-h}e_iv^{h}=v^{\alpha_i(h)}e_i,v^{-h}f_iv^{h}=v^{-\alpha_i(h)}f_i$ for any $i\in I$ and $h\in\MP$,
\item $e_if_j-f_je_i=\delta_{ij}(K_i-K_i^{-1})/(v_i-v_i^{-1})$ for any $i,j\in I$,
\item\label{qSerre1} $\sum_{k=0}^{1-a_{ij}}(-1)^k e_i^{(k)}e_je_i^{(1-a_{ij}-k)}=0$
for any $i,j\in I$ with $i\ne j$,
\item\label{qSerre2} $\sum_{k=0}^{1-a_{ij}}(-1)^k f_i^{(k)}f_jf_i^{(1-a_{ij}-k)}=0$
for any $i,j\in I$ with $i\ne j$,
\ee
where $K_i=v^{d_ih_i}, v_i=v^{d_i}$ and $e_i^{(n)}=e_i^n/[n]_{d_i}!,f_i^{(n)}=f_i^n/[n]_{d_i}!$.
\edf

\Rem
$U_v$ has $\cor$-antiinvolution $\sigma$ and $\Q$-antiinvolution $\Omega,\Upsilon$  defined by
\begin{align*}
\sigma(e_i) &= f_i,\quad
\sigma(f_i)=e_i,\quad
\sigma(v^h)=v^h,\\
\Omega(e_i) &= f_i,\quad
\Omega(f_i)=e_i,\quad
\Omega(v^h)=v^{-h},\quad
\Omega(v)=v^{-1},\\
\Upsilon(e_i) &= v_if_iK_i^{-1},\quad
\Upsilon(f_i)=v_i^{-1}K_ie_i,\quad
\Upsilon(v^h)=v^{-h},\quad
\Upsilon(v)=v^{-1}.
\end{align*}
\enrem

\Prop[{\cite[\S3.2]{Lus}}]\label{QPBW}
Define the three $\cor$-subalgebras $U_v^\pm,U_v^{0}$ of $U_v$ by
\begin{align*}
U_v^+=\langle e_i\mid i\in I\rangle,\quad
U_v^-=\langle f_i\mid i\in I\rangle,\quad
U_v^0=\langle v^h\mid h\in\MP\rangle.
\end{align*}
\bna
\item the canonical map
$U^{-}_v\otimes_{\cor}U^{0}_v\otimes_{\cor}U^{+}_v\to U_v$
is a $\cor$-vector space isomorphism,
\item
$U^{0}_v$ is canonically isomorphic to the group $\cor$-algebra $\cor[\MP]$.
\ee
\enprop

\Def
We define the following two maps both are assured by Proposition \ref{QPBW}.
\bna
\item the Harish-Chandra projection $\HC:U_v\to U_v^0$ which is a $\cor$-linear projector
from $U_v=U_v^0\oplus((\sum_{i\in I}f_iU_v)+(\sum_{i\in I}U_ve_i))$ to $U_v^0$.
\item the evaluation map $\EV_{\lambda}:U_v^0\to \cor$ which is 
an $\cor$-algebra homomorphism determined by the assignment $\EV_{\lambda}(v^h)=v^{\lambda(h)}$ for each $h\in\MP$.
\ee
\edf

\subsection{Shapovalov forms}
\Def[{\cite[\S3.4.5]{Lus}}]
For $\lambda\in\MP$, the Verma module is a left $U_v$-module 
\begin{align*}
M(\lambda) = U_v/(\sum_{h\in\MP}U_v(v^h-v^{\lambda(h)})+\sum_{i\in I}U_ve_i).
\end{align*}
\edf

\Rem
In virtue of Proposition \ref{QPBW}, as a $U^-_v$-module $M(\lambda)$ is free of rank 1 and we see that
$M(\lambda)$ has a proper maximum $U_v$-submodule $K(\lambda)$.
We reserve the symbol $\varphi_\lambda$ for the $U^-_v$-linear isomorphism 
$\varphi_\lambda:U^-_v\stackrel{\sim}{\to} M(\lambda),
u\mapsto u+\sum_{h\in\MP}U_v(v^h-v^{\lambda(h)})+\sum_{i\in I}U_ve_i$.
\enrem

\Def
For $\lambda\in\MP$, we define an irreducible $U_v$-module $V(\lambda)=M(\lambda)/K(\lambda)$.
\edf

\Prop\label{exi_sh}
Let us take $\lambda\in\MP$.
\bna
\item There exists a unique non-degenerate symmetric $\cor$-bilinear form $\langle,\rangle_{\SH}:V(\lambda)\times V(\lambda)\to \cor$ 
(called the Shapovalov form)
such that
$\langle 1_{\lambda},1_{\lambda}\rangle_{\SH}=1$ and
$\langle uv,w\rangle_{\SH}=\langle v,\sigma(u)w\rangle_{\SH}$ for all $v,w\in V(\lambda)$ and $u\in U_v$,
\item There exist unique bi-additive non-degenerate maps $\langle,\rangle_{\QSH}:V(\lambda)\times V(\lambda)\to \cor$ and
$\langle,\rangle_{\RSH}:V(\lambda)\times V(\lambda)\to \cor$ such that for all $X\in\{\QSH,\RSH\}$
\bnum
\item $\langle aw_1,w_2\rangle_{X}=\sigma(a)\langle w_1,w_2\rangle_{X}$,
$\langle w_1,aw_2\rangle_{X}=a\langle w_1,w_2\rangle_{X}$ and 
$\langle w_1,w_2\rangle_{X}=\sigma(\langle w_2,w_1\rangle_{X})$,
\item $\langle 1_{\lambda},1_{\lambda}\rangle_{X}=1$ and 
$\langle uw_1,w_2\rangle_{\QSH}=\langle w_1,\Omega(u)w_2\rangle_{\QSH}$, $\langle uw_1,w_2\rangle_{\QSH}=\langle w_1,\Upsilon(u)w_2\rangle_{\RSH}$.
\ee
for all $w_1,w_2\in V(\lambda)$, $u\in U_v$ and $a\in\cor$.
\ee
\enprop

\Proof
All cases are similar and standard. % (for Shapovalov forms, see hoge). 
Let $X\in\{\SH,\QSH,\RSH\}$ and put $\Xi=\sigma,\Omega,\Upsilon$ according to $X=\SH,\QSH,\RSH$ respectively.
We just sketch the constructions of $\langle,\rangle_{X}$ since they are needed 
in the proof of Proposition \ref{SHANDQSH}.

Define $\langle,\rangle'_{X}:M(\lambda)\times M(\lambda)\to \cor$ 
by $\langle w_1,w_2\rangle'_{X}=\EV_{\lambda}(\HC(\Xi(\varphi^{-1}_{\lambda}(w_1))\varphi^{-1}_{\lambda}(w_2)))$. 
Then, we see that $\RAD(\langle,\rangle'_{X})=K(\lambda)$ and
%Now it is clear that 
$\langle,\rangle'_{X}$ decent to 
$\langle,\rangle_{X}:V(\lambda)\times V(\lambda)\to\cor$ which satisfy the desired conditions.
\QED

\subsection{Lusztig lattices}
\Th[{\cite[Theorem 14.4.3]{Lus}},{\cite[Theorem 4.5]{Lus2}}]\label{LUS2}
Let $U_v^{\MA}$ be an $\MA$-subalgebra generated by
$\{e_i^{(n)},f_i^{(n)},K_i^{\pm1}\mid i\in I,n\geq 0\}$.
Then, $U_v^{\MA}$ is an $\MA$-lattice of $U_v$.
\enth

%\Rem
%$U^{\MA}$ has a triangular decomposition compatible with Proposition \ref{QPBW}.
%In particular, $U^{+,\MA}=\langle e^{(r)}_i\mid i\in I,r\geq 0\rangle_{\textrm{$\MA$-alg}}$ and
%$U^{-,\MA}=\langle f^{(r)}_i\mid i\in I,r\geq 0\rangle_{\textrm{$\MA$-alg}}$
%where $U^{\pm,\MA}=U^{\MA}\cap U_v^{\pm}$.
%\enrem

%More precisely, let $U_\MA^{+}$ (resp. $U_\MA^{-}$) be an $\MA$-subalgebra of $U^{+}_v$ (resp. $U^{-}_v$) 
%generated by the divided powers $\{e_i^{(n)}\mid i\in I, n\geq 0\}$ (resp. $\{f_i^{(n)}\mid i\in I, n\geq 0\}$) and
%let $U_{\MA}^{0}$ be an $\MA$-subalgebra of $U^{0}_v$ generated by $K_i^{\pm1}$

%Recall that a $U_v$-module $M$ is called integrable (see ~\cite[\S3.5]{Lus}) if
\Def[{\cite[\S3.5]{Lus}}]
A $U_v$-module $M$ is called integrable 
\bnum
\item $M=\oplus_{\nu\in\Gh^\ast}M_{\nu}$ with $\dim M_{\nu}<+\infty$ where
$M_{\nu} = \{ m\in M\mid \forall h\in\MPC,v^hm=v^{\nu(h)}m \}$,
\item for any $m\in M$ and $i\in I$, there exists some $n>0$ such that $f_i^nm=e_i^nm=0$.
%we have $\dim_{\Q(v_i)} U_{v,i}m<+\infty$ where $U_{v,i}$ is a $\Q(v_i)$-subalgebra of $U_v$
%generated by $\{e_i,f_i,K_i^{\pm1}\}$.
\ee
\edf

\Rem[{\cite[Proposition 3.5.8, Proposition 5.2.7]{Lus}}]
It is well-known that %(see ~\cite[Proposition 3.5.8, Proposition 5.2.7]{Lus})
\bna
\item for $\lambda\in\MP$, $V(\lambda)$ is integrable if and only if $\lambda\in\MP^+$.
\item for $\lambda\in\MP^+$, the set of weights $P(\lambda):=\{\mu\in\Gh^\ast\mid V(\lambda)_{\mu}\ne 0\}$ of $V(\lambda)$
is $W$-invariant.
\ee
\enrem

\Th[{\cite[Theorem 14.4.11]{Lus}}]\label{LUS} Assume $\lambda\in\MP^+$. 
Then, $V(\lambda)^{\MA}:=U_v^{\MA}1_{\lambda}$ is an $\MA$-lattice of $V(\lambda)$ compatible with the
weight space decomposition of $V(\lambda)$. In
other words, $V(\lambda)_{\nu}^{\MA}:=V(\lambda)_{\nu}\cap V(\lambda)^{\MA}$ is an $\MA$-lattice of $V(\lambda)_{\nu}$
for all $\nu\in P(\lambda)$ and we have $V(\lambda)^{\MA}=\oplus_{\nu\in P(\lambda)}V(\lambda)_{\nu}^{\MA}$. 
\enth

\Def
For $\lambda\in \MP^+$ and $\mu\in P(\lambda)$, we define $\SH_{\lambda,\mu} = \det\SHM_{\lambda,\mu}\in \cor^\times/v^{2\Z}$ and
\begin{align*}
\SHM_{\lambda,\mu} &= \SHM_{\lambda,\mu}(A) = \GRM_{\cor,\ID_{\cor},\MA}(V(\lambda)_{\mu},\SH,V(\lambda)^{\MA}_{\mu}),\\
\QSHM_{\lambda,\mu} &= \QSHM_{\lambda,\mu}(A) = \GRM_{\cor,\BAR,\MA}(V(\lambda)_{\mu},\QSH,V(\lambda)^{\MA}_{\mu}),\\
\RSHM_{\lambda,\mu} &= \RSHM_{\lambda,\mu}(A) = \GRM_{\cor,\BAR,\MA}(V(\lambda)_{\mu},\RSH,V(\lambda)^{\MA}_{\mu}).
\end{align*}
\edf

\Rem\label{speci}
Theorem \ref{LUS2} and Theorem \ref{LUS} allow us to specialize the integrable highest weight modules.
For any $\xi\in\C^{\times}$ and $\lambda\in\MP^+$, the specialized module $V(\lambda)|_{v=\xi}:=\C_{\xi}\otimes_{\MA} V(\lambda)^{\MA}$
is a module over $U_v|_{v=\xi}:=\C_{\xi}\otimes_{\MA}U_v^{\MA}$ where $\C_{\xi}=\C$ is an $\MA$-module where $v$ acts as the multiplication of $\xi$.
\enrem

\Rem\label{speci2}
For $\lambda\in \MP^+$ and $\xi\in \C^{\times}$,
$V(\lambda)|_{v=\xi}$ is irreducible if and only if $\SH_{\lambda,\mu}(\xi)\ne 0$ for all $\mu\in P(\lambda)$.
Thus, for a generic $\xi\in \C^{\times}$, $V(\lambda)|_{v=\xi}$ is irreducible.
%As far as the author knows, this fact was first pointed out in ~\cite[Theorem 4.12]{Lus3} by a different method.
It is expected that $V(\lambda)|_{v=\xi}$ is irreducible when $\xi$ is not a root of unity~\cite[Problem 2]{Kas}.
\enrem

\Prop\label{SHANDQSH}
Let $\lambda\in\MP^+$ and take $\mu\in P(\lambda)$. We may assume $\SHM_{\lambda,\mu}=\QSHM_{\lambda,\mu}$ and
$D\QSH^M_{\lambda,\mu}=\RSH^M_{\lambda,\mu}$ for some diagonal matrix $D$ all of whose diagonal entries belong to $v^{\Z}$.
\enprop

\Proof
By the constructions mentioned in the proof of Proposition \ref{exi_sh},
it is enough to show that we can choose an $\MA$-basis of $V(\lambda)^{\MA}_\mu$ as a subset of %consisting of the elements in
$\{f^{(n_1)}_{i_1}\cdots f^{(n_s)}_{i_s}1_{\lambda}\mid\substack{s\geq 0\\ 1\leq \forall j\leq s,n_j\geq 1,i_j\in I}\}$.
In virtue of ~\cite[Theorem 6.5]{Lak}, it is possible\footnote{For the basic representations 
of $U_v(A^{(1)}_{n}), U_v(A^{(2)}_{2n-1}), U_v(A^{(2)}_{2n}), U_v(D^{(1)}_{n}), U_v(D^{(2)}_{n+1}), U_v(B^{(1)}_{n})$,
we can choose such monomial $\MA$-bases using combinatorial algorithms~\cite{LLT,LT,KKw}.}.

\QED

\Cor\label{detichi}
For $\lambda\in\MP^+$ and  $\mu\in P(\lambda)$, we have $\det\QSH^M_{\lambda,\mu}=\det\RSH^M_{\lambda,\mu}$ in $\cor^{\times}$.
\encor

\Proof
We may assume $D\QSH^M_{\lambda,\mu}=\RSH^M_{\lambda,\mu}$ for some $D$ with $\det D\in v^{\Z}$.
Since $\TRANS{\QSH^M_{\lambda,\mu}}=\BAR(\QSH^M_{\lambda,\mu})$ and $\TRANS{\RSH^M_{\lambda,\mu}}=\BAR(\RSH^M_{\lambda,\mu})$,
we easily see that $\det D=\BAR(\det D)$.
\QED

\Prop\label{Weyl}
For $\lambda\in\MP^+$, $i\in I$ and $\mu\in P(\lambda)$, we may assume $\QSHM_{\lambda,\mu}=\QSHM_{\lambda,s_i(\mu)}$.
\enprop

\Proof
Since $V(\lambda)$ is integrable, % when $\lambda$ is dominant integral,
we have for each $i\in I$
mutually inverse linear operators $T'_{i,-1}=\sum_{\nu\in P(\lambda)}T'^{(\nu)}_{i,-1}$ and $T''_{i,+1}=\sum_{\nu\in P(\lambda)}T''^{(\nu)}_{i,+1}$ 
on $V(\lambda)$ where (see ~\cite[\S5.2]{Lus})
%where $T_{i,\nu}$ and ${}^{\omega}T'_{i,\nu}$ are Lusztig's braid group action defined by
\begin{align*}
T'^{(\nu)}_{i,-1}: V(\lambda)_{\nu}\isoto V(\lambda)_{s_i(\nu)}, \quad w\MAPSTO \sum_{\substack{a,b,c\geq 0 \\ a-b+c=\nu(h_i)}}(-1)^bv_i^{ac-b}f_i^{(a)}e_i^{(b)}f_i^{(c)}w, \\
T''^{(\nu)}_{i,+1}: V(\lambda)_{\nu}\isoto V(\lambda)_{s_i(\nu)}, \quad w\MAPSTO \sum_{\substack{a,b,c\geq 0 \\ -a+b-c=\nu(h_i)}}(-1)^bv_i^{b-ac}e_i^{(a)}f_i^{(b)}e_i^{(c)}w.
\end{align*}

It is clear that we have $\langle T'_{i,-1}(v),T'_{i,-1}(w)\rangle_{\QSH}=\langle v,T''_{i,+1}(T'_{i,-1}(w))\rangle_{\QSH}=\langle v,w\rangle_{\QSH}$.
Since $T'_{i,-1}$ (and $T''_{i,+1}$) preserves $V(\lambda)^\MA$~\cite[Proposition 41.2.2]{Lus}, we are done.
\QED

%\Cor\label{Weyl}
%For $\lambda\in\MP^+$ and $\mu_1,\mu_2\in P(\lambda)$, we have $\SH_{\lambda,\mu_1}=\SH_{\lambda,\mu_2}$ if $W\mu_1=W\mu_2$.
%\encor

\section{Shapovalov determinants of the basic representations}% of affine quantum groups}
\subsection{Untwisted affine A,D,E case}

\begin{figure}
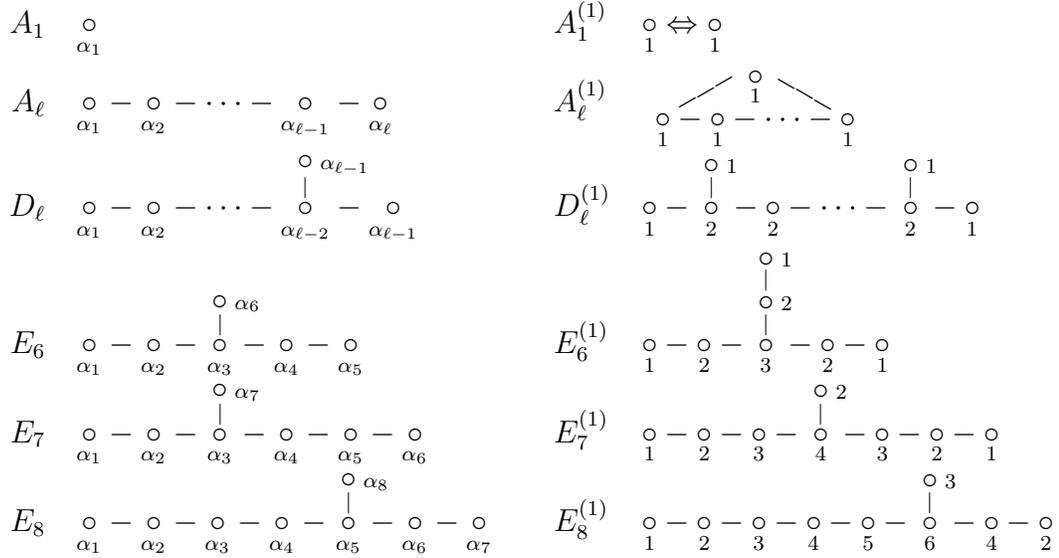

\[
\begin{array}{r@{\quad}l@{\qquad}l@{\quad}l}
A_1&\node{}{\alpha_1}&
A_1^{(1)} & \node{}{1} \Leftrightarrow \node{}{1} \\
A_\ell & \node{}{\alpha_1} - \node{}{\alpha_2} - \cdots - \node{}{\alpha_{\ell-1}}-\node{}{\alpha_\ell} &
A_\ell^{(1)}  & 
\begin{array}{c}
\raisebox{-12pt}{\rotatebox{30}{$-\!\!-\!\!-$}}\node{}{1}\raisebox{0pt}{\rotatebox{-30}{$-\!\!-\!\!-$}} \\[-7pt]
\node{}{1}-\node{}{1}-\cdots-\node{}{1} 
\end{array}\\
D_\ell& \node{}{\alpha_1} - \node{}{\alpha_2} - \cdots - \node{\ver{}{\alpha_{\ell-1}}}{\alpha_{\ell-2}} -\node{}{\alpha_{\ell-1}} &
D_\ell^{(1)}  & \node{}{1}-\node{\ver{}{1}}{2}-\node{}{2}-\cdots-\node{\ver{}{1}}{2}-\node{}{1} \\
E_6& \node{}{\alpha_1} - \node{}{\alpha_2} - \node{\ver{}{\alpha_{6}}}{\alpha_{3}} -\node{}{\alpha_{4}} -\node{}{\alpha_5} &
E_6^{(1)}  & \node{}{1}-\node{}{2}-\node{\overset{\ver{}{1}}{\ver{}{2}}}{3}-\node{}{2}-\node{}{1} \\
E_7& \node{}{\alpha_1} - \node{}{\alpha_2} - \node{\ver{}{\alpha_{7}}}{\alpha_{3}} -\node{}{\alpha_{4}} -\node{}{\alpha_5} -\node{}{\alpha_6} &
E_7^{(1)} & \node{}{1}-\node{}{2}-\node{}{3}-\node{\ver{}{2}}{4}-\node{}{3}-\node{}{2}-\node{}{1}\\
E_8& \node{}{\alpha_1} - \node{}{\alpha_2} - \node{}{\alpha_{3}} -\node{}{\alpha_{4}} -\node{\ver{}{\alpha_{8}}}{\alpha_5} -\node{}{\alpha_6}-\node{}{\alpha_7} &
E_8^{(1)} & \node{}{1}-\node{}{2}-\node{}{3}-\node{}{4}-\node{}{5}-\node{\ver{}{3}}{6}-\node{}{4}-\node{}{2} \\
\end{array}
\]
\caption{Finite and untwisted affine Dynkin diagrams of A,D,E.}
\label{untwisted}
\end{figure}

Let $X=(a_{ij})_{i,j\in I}$ be a Cartan matrix of type A,D,E %of rank $\ell=|I|$ 
and let $\HA=X^{(1)}$ be the extended Cartan matrix of $X$ indexed by $\HI=\{0\}\sqcup I$ as in Figure \ref{untwisted}.
Let $(a_i)_{i\in\HI}$ be the numerical labels of $\HA$ in Figure ~\ref{untwisted} and let $\delta=\sum_{i\in \HI}a_i\alpha_i$.

%The purpose of this section is to review two facts 
%on the $U_v(\HA)$-module $V(\Lambda_0)$ (also known as the basic module) 
%which will be used in our proof of Theorem ~\ref{MainTheorem}.

\Rem[{\cite[\S12.6]{Kac}}]
Let $W=W(\HA)$. We have
\begin{align}\label{Wei}
P(\Lambda_0)=\{w\Lambda_0-d\delta\mid w\in W, d\geq 0\}
(=\{w\Lambda_0-d\delta\mid w\in W/W_0, d\geq 0\})
\end{align}
Note that $W/W_0\cong W(X)$ where $W_0=\{w\in W\mid w\Lambda_0=\Lambda_0\}$.
\enrem

We shall very briefly review the vertex operator construction of $V(\Lambda_0)$, i.e., an
explicit realization of $V(\Lambda_0)$
as the tensor product of a polynomial algebra and a (twisted) group algebra~\cite{FJ} (however, we will mainly refer ~\cite{CJ} instead).
Recall the Drinfeld new realization of $U_v(\HA)$
which gives loop-like generators~\cite{Dri}.% of $U_v(\HA)$.

\Th\label{DNEW}
As a $\cor$-algebra,
$U_v(\HA)$ is isomorphic to the $\cor$-algebra generated by
$\{x^{\pm}_{i,s}\mid i\in I,s\in\Z\}\cup\{h_{i,r}\mid i\in I,r\in\Z\setminus\{0\}\}\cup \{C^{\pm1},D^{\pm1},K^{\pm1}_i\mid i\in I\}$ with
the following defining relations.
\bna
\item $C^{\pm1}$ are central and $CC^{-1}=C^{-1}C=DD^{-1}=D^{-1}D=K_iK_i^{-1}=K_i^{-1}K_i=1$,
\item $K_iK_j=K_jK_i$, $K_ih_{j,r}=h_{j,r}K_i$, $K_ix^{\pm}_{j,s}K^{-1}_i=v^{\pm a_{ij}}x^{\pm}_{j,s}$,
\item $DK_i=K_iD$, $Dh_{j,r}D^{-1}=v^rh_{j,r}$, $Dx^{\pm}_{j,s}D^{-1}=v^rx^{\pm}_{j,s}$,
\item $\displaystyle[h_{i,r},h_{j,r'}]=\delta_{r+r',0}\frac{[ra_{ij}]}{r}\frac{C^r-C^{-r}}{v-v^{-1}}$,
$\displaystyle[x^+_{i,s},x^-_{j,s'}]=\delta_{i,j}\frac{C^{-s'}\psi^+_{i,s+s'}-C^{-s}\psi^-_{i,s+s'}}{v-v^{-1}}$,
\item $\displaystyle[h_{i,\pm r},x^{\pm}_{j,s}]=\pm \frac{[ra_{ij}]}{r}x^{\pm}_{j,s\pm r}$, 
$\displaystyle[h_{i,\mp r},x^{\pm}_{j,s}]=\pm \frac{[ra_{ij}]}{r}C^rx^{\pm}_{j,s\pm r}$,
\item $\displaystyle x^{\pm}_{i,s+1}x^{\pm}_{j,s'}-v^{\pm a_{ij}}x^{\pm}_{j,s'}x^{\pm}_{i,s+1}=v^{\pm a_{ij}}x^{\pm}_{i,s}x^{\pm}_{j,s'+1}-x^{\pm}_{j,s'+1}x^{\pm}_{i,s}$,
\item $\sum_{\pi\in\Sym_{1-a_{pq}}}\sum_{\substack{\BR=(r_{\ell})\in\Z^{1-a_{pq}} \\ 0\leq k\leq 1-a_{pq}}}(-1)^k\Bino{1-a_{pq}}{k}x^{\pm}_{p,r_{\pi(1)}}\cdots x^{\pm}_{p,r_{\pi(k)}}x^{\pm}_{q,s}x^{\pm}_{p,r_{\pi(k+1)}}\cdots x^{\pm}_{p,r_{\pi(m)}}=0$,
\ee
for all $i,j,p,q\in I$ with $p\ne q$ and $r,r'\in\Z\setminus\{0\}, s,s'\in\Z$
where we set $\psi^+_{-s}=\psi^-_{s}=0$ for $s\geq 1$ and
%\begin{align*}
$\sum_{t\geq 0}\psi^\pm_{i,\pm t}z^{\pm t}=K^{\pm 1}_i\exp(\pm (v-v^{-1})\sum_{r\geq 1}h_{i,\pm r}z^{\pm r})$.
%\end{align*}
\enth

\Rem\label{Dchoice}
Actually, as shown by Beck~\cite{Bec} (see also ~\cite[Theorem 1]{BCP}), we may choose the generators in Theorem ~\ref{DNEW}
as particular elements in $U_v(\HA)$, some of them are constructed utilizing Lusztig's braid group action on $U_v(\HA)$
from the Chevalley generators of $U_v(\HA)$.
Although we omit the construction, in the rest of this paper we assume this special choice for each generator in Theorem ~\ref{DNEW}.
The choice is used essentially in that
\bna
\item an $\MA$-subalgebra generated by $\{(x_{i,n}^{\pm})^r/[r]!\mid i\in I, r\geq 0\}$ coincides with $U_v^{\pm,\MA}$~\cite[Corollary 2.2]{BCP}
where $U_v^{+,\MA}=\langle e^{(r)}_i\mid i\in I,r\geq 0\rangle_{\textrm{$\MA$-alg}}$ and
$U_v^{-,\MA}=\langle f^{(r)}_i\mid i\in I,r\geq 0\rangle_{\textrm{$\MA$-alg}}$.
This is used in the proof of Theorem ~\ref{VCON} (ii),
\item we have $\Omega(h_{i,\pm r})=h_{i,\mp r}$ for $i\in I$ and $r\geq 1$~\cite[Lemma 3]{BCP}. This will be used in the proof of Theorem ~\ref{MainTheorem}.
\ee
\enrem

\Rem
Let $U^{-}(0)$ be the $\cor$-subalgebra of $U_v(\HA)$ generated by $\{h_{i,-r}\mid i\in I, r\geq 1\}$.
Clearly, it is a commutative algebra.
In virtue of ~\cite[Proposition 1.3]{BCP}, $\{h_{i,-r}\mid i\in I, r\geq 1\}$ is algebraically independent over $\cor$.
\enrem
%, i.e.,$U^{-}(0)=\cor[h_{i,-r}\mid i\in I, r\geq 1]$ is a polynomial ring of variables $\{h_{i,-r}\mid i\in I, r\geq 1\}$ over $\cor$.
%The following is a (partial) realization of $V(\Lambda_0)$ and its Lusztig lattice $V(\Lambda_0)^\MA$.

\Th[{\cite[Theorem2, Theorem 3]{CJ}}]\label{VCON}
Let $Q=\oplus_{p\in I}\Z\alpha_p$ be the root lattice of $X$. %and regard it as an additive subgroup of the weight lattice $\MP$ of $\HA$.
We can put a $U_v(\HA)$-module structure on a $\cor$-vector space
$V:=U^-(0)\otimes_{\cor}\cor[Q]$ %=\cor[h_{i,-r}\mid i\in I, r\geq 1]\otimes_{\cor}\cor[Q]$ where $Q:=\oplus_{p\in I}\Z\alpha_p$ 
%into a $U_v(\HA)$-module 
by extending the following assignment
of action $\ACT$ on $V$
\begin{align*}
h_{i,-s}\ACT(h_{i_1,-r_1}\cdots h_{i_j,-r_j}\otimes e^{\beta}) &= h_{i,-s}h_{i_1,-r_1}\cdots h_{i_j,-r_j}\otimes e^{\beta}, \\
h_{i,s}\ACT(h_{i_1,-r_1}\cdots h_{i_j,-r_j}\otimes e^{\beta}) &= \sum_{k=1}^{j} h_{i_1,-r_1}\cdots h_{i_{k-1},-r_{k-1}}\frac{\delta_{s,r_k}[sa_{i,i_{k}}][s]}{s}h_{i_{k+1},-r_{k+1}}\cdots h_{i_j,-r_j}\otimes e^{\beta},
\end{align*} 
to all the generators in Theorem \ref{DNEW} and can obtain a $U_v(\HA)$-module isomorphism $\varphi:V\isoto V(\Lambda_0)$ such that
\bnum
\item $\varphi(h_{i_1,-r_1}\cdots h_{i_j,-r_j}\otimes e^{\beta})\in V(\Lambda_0)_{\Lambda_0-\beta-(\sum_{k=1}^{j}r_k+(\beta,\beta)/2)\delta}$,
\item $\varphi(\MA[\MMP_{i,r}\mid i\in I,r\geq 1]\otimes_{\MA}\MA[Q])=V(\Lambda_0)^\MA$
where $1+\sum_{r\geq 1}\MMP_{i,r}z^r=\exp(\sum_{r\geq 1}\frac{h_{i,-r}}{[r]}z^r)$.
\ee
Here $(i_1,\cdots,i_j)\in I^j,(r_1,\cdots,r_j)\in \NNN^j,\beta=\sum_{p\in I}b_p\alpha_p\in Q$, $i\in I$, $s\geq 1$,
$(\beta,\beta)=\sum_{p,q\in I}b_pa_{pq}b_q$ and
we regard $Q$ as an additive subgroup of the weight lattice $\MP$ of $\HA$.
\enth

\Th\label{MainTheorem}
Let $X=(a_{ij})_{i,j\in I}$ be a Cartan matrix of type A,D,E and
let $\HA$ be the extended Cartan matrix of $X$ indexed by $\HI=\{0\}\sqcup I$ as in Figure ~\ref{untwisted}.
We have 
\begin{align*}
\SH_{\Lambda_0,w\Lambda_0-d\delta}(\HA) = 
\prod_{s=1}^{d}(\det [X]_s)^{\sum_{\lambda\in\PAR(d)}\frac{m_s(\lambda)}{|I|}N_{\lambda}}
\end{align*}
in $\cor^\times/v^{2\Z}$ for all $w\in W$ and $d\geq 0$
where $N_{\lambda}=\prod_{u\geq 1}\binom{m_u(\lambda)+|I|-1}{m_u(\lambda)}$ and $[X]_s=([a_{ij}]_s)_{i,j\in I}$.
\enth

\Proof
By Proposition \ref{SHANDQSH} and Proposition \ref{Weyl}, 
we shall calculate $\det\QSHM_{\Lambda_0,\Lambda_0-d\delta}(\HA)$ instead of $\SH_{\Lambda_0,w\Lambda_0-d\delta}(\HA)$.
By Theorem \ref{VCON}, we see
\bnum
\item $V(\Lambda_0)_{\Lambda_0-d\delta}$ has a $\cor$-basis $\{
h'_{i_1,-r_1}\cdots h'_{i_j,-r_j}\otimes e^{0}\mid 
\textrm{$(i_\ell)\in I^j,(r_\ell)\in \NNN^j$ with $\sum_{k=1}^{j}r_k=d$}
%\sum_{k=1}^{j}r_k=d
\}$,
\item $V(\Lambda_0)^{\MA}_{\Lambda_0-d\delta}$ has an $\MA$-basis $\{
\MMP_{i_1,r_1}\cdots \MMP_{i_j,r_j}\otimes e^{0}\mid 
\textrm{$(i_\ell)\in I^j,(r_\ell)\in \NNN^j$ with $\sum_{k=1}^{j}r_k=d$}
%\sum_{k=1}^{j}r_k=d
\}$,
\ee
where $h'_{i,\pm r}=h_{i,\pm r}/[r]$ for $i\in I$ and $r\geq 1$. 
Note that $1+\sum_{r\geq 1}\MMP_{i,r}z^r=\exp(\sum_{r\geq 1}h'_{i,-r}z^r)$.

Since $\Omega(h'_{i,-r})=h'_{i,r}$ (see Remark \ref{Dchoice}) and 
\begin{align*}
(sh'_{i,s})\ACT(h'_{i_1,-r_1}\cdots h'_{i_j,-r_j}\otimes e^{0}) = 
\sum_{k=1}^{j} \delta_{s,r_k}[a_{i,i_{k}}]_sh'_{i_1,-r_1}\cdots h'_{i_{k-1},-r_{k-1}}h'_{i_{k+1},-r_{k+1}}\cdots h'_{i_j,-r_j}\otimes e^{0},
\end{align*}
we can apply Corollary \ref{WITHOUTK} and get
$\det\QSHM_{\Lambda_0,\Lambda_0-d\delta}(A) = \prod_{s=1}^{d}(\det [A]_s)^{\sum_{\lambda\in\PAR(d)}\frac{m_s(\lambda)}{|I|}N_{\lambda}}$.
\QED

\Cor\label{CJcomp}
Let $\Ggz$ be a finite-dimensional simple simply-laced complex Lie algebra. 
\bna
\item
If $\xi\in\C^{\times}$ is not a root of unity, then $V(\Lambda_0)|_{v=\xi}$ is 
irreducible.
\item Let $\Ggz=\Gg(X)$ and let $\ell\geq 1$. $V(\Lambda_0)|_{v=\exp(\frac{2\pi\sqrt{-1}}{\ell})}$ is 
irreducible if and only if %where $\xi=\exp(\frac{2\pi\sqrt{-1}}{\ell})$ if and only if 
\bnum
\item $\GCD(\ell,2n)=1,2$ when $X=A_{n-1}$ for $n\geq 2$,
%\item $\GCD(\ell,n')=1$ and $\ell\not\in 4\Z$ when $X=A_{n-1}$ for even $n=2^{e}n'\geq 2$ where $e\geq 1$ and $n'\not\in 2\Z$,
\item $\ell\not\in 4\Z$ (resp. $\not\in 3\Z,\not\in 4\Z,\not\in 60\Z$) when $X=D_{n}$ for $n\geq 4$ (resp. $E_6,E_7,E_8$).
\ee
\ee
\encor

\Proof
The specialized module is irreducible if and only if $(\det [X]_k)|_{v=\xi}\ne 0$ for all $k\geq 1$.
Type by type, $\det ([X]_1)$ is given by $\det ([A_{n-1}]_1)=[n]$, $\det ([D_m]_1)=v^{-m}\Phi_{4}(v)\Phi_{4}(v^{m-1})$, 
$\det ([E_6]_1) = v^{-6}\Phi_{3}(v^2)\Phi_{24}(v)$, $\det ([E_7]_1) = v^{-7}\Phi_{4}(v)\Phi_{36}(v)$,
$\det ([E_8]_1) = v^{-8}\Phi_{60}(v)$
for $n\geq 2$, $m\geq 4$ where $\Phi_{j}(x)\in\Z[x]$ is the $j$-th cyclotomic polynomial.
\QED

%\Rem
%Chari-Jing showed that if the Coxeter number $N$ of $\Ggz$ is coprime to $\ellgeq 1$, then 
%$V(\Lambda_0)|_{v=\exp(\frac{2\pi\sqrt{-1}}{\ell})}$ is irreducible~\cite[Theorem 4]{CJ}.
%Recall that $N=n,2(m-1),12,18,30$ when $X=A_{n-1}$,$D_m$, $E_6$, $E_7$, $E_8$ respectively for $n\geq 2$, $m\geq 4$.
%However, as seen in the proof of ~\cite[Theorem 4]{CJ} they actually prove if part of Corollary \ref{CJcomp} because
%what their proof requires is only that $(\det [X]_k)|_{v=\xi}\ne 0$ for all $k\geq 1$
%\enrem

\Rem
Let $\cor$ be a field of characteristic $p\geq 2$.
DeConcini-Kac-Kazhdan proved that
the modulo-$p$ reduction $\cor\otimes_{\Z}V(\Lambda_0)^{\Z}$ of the basic $\Ggg$-module $V(\Lambda_0)$ where
$V(\Lambda_0)^\Z$ is the Kostant $\Z$-form of $V(\Lambda_0)$ remains irreducible 
if and only if $\det X\ne 0$ in $\cor$~\cite[Corollary 3.1]{DcKK}.
Corollary \ref{CJcomp} can be seen a quantum analog of
DeConcini-Kac-Kazhdan's result.
Note that when $X=E_8$ the modulo-$p$ reduction preserves irreduciblity for any prime $p\geq 2$ because $\det E_8=1$.
Since the quantum characteristic $\min\{k\geq 1\mid [k]|_{v=\xi}=0\}$ takes arbitrary value $\geq 2$,
quantum case is more subtle than the classical case.% treated by DeConcini-Kac-Kazhdan.
\enrem

\subsection{Twisted affine A,D,E case}

Let $\HA=(a_{ij})_{i,j\in\HI=\{\varepsilon\}\sqcup I}$ be a twisted affine {\GCM} as in Figure ~\ref{twisted}
%where $\varepsilon=0$ if $\HA\ne A^{(2)}_{2\ell}$ for any $\ell\geq 1$ and otherwise $\varepsilon=\ell$.
where $\varepsilon=\ell$ if $\HA=A^{(2)}_{2\ell}$ and otherwise $\varepsilon=0$.
Let $\delta=\sum_{i\in \HI}a_i\alpha_i$ where $(a_i)_{i\in\HI}$ are the numerical labels of $\HA$.
As in untwisted cases (see Remark \ref{Wei}),
it is known that $P(\Lambda_{\varepsilon})=\{w\Lambda_{\varepsilon}-d\delta\mid w\in W(\HA), d\geq 0\}$~\cite[\S12.6]{Kac}.

We discuss on $\SH_{\Lambda_{\varepsilon},w\Lambda_{\varepsilon}-d\delta}(\HA)$ and
give a conjectural formula for it. %$\SH_{\Lambda_0,w\Lambda_0-d\delta}(\HA)$.
Our motivation comes from the fact that %when $p$ is an odd prime number
the Cartan determinant of a faithful $p$-block of the Schur double covers of the symmetric groups 
coincides with $\SH_{\Lambda_{\ell},w\Lambda_{\ell}-d\delta}(A^{(2)}_{p-1})|_{v=1}$ for suitable $w\in W(A^{(2)}_{p-1})$ and $d\geq 0$
when $p=2\ell+1$ is an odd prime~\cite[Theorem 5.6]{BO2}.

\begin{figure}
\[
\begin{array}{r@{\quad}l@{\qquad}l@{\quad}l}
A_2^{(2)} &\node{1}{\alpha_0} {\quad\!\!\!\!\!\!\!\!}\Rrrightarrow{\quad\!\!\!\!\!} \node{2}{\alpha_1} &
D_{\ell+1}^{(2)}  &\node{1}{\alpha_0}\Leftarrow \node{1}{\alpha_1}-\cdots-\node{1}{\alpha_{\ell-1}}\Rightarrow\node{1}{\alpha_\ell} \\
A_{2\ell}^{(2)}  &\node{1}{\alpha_0}\Rightarrow \node{2}{\alpha_1}-\cdots-\node{2}{\alpha_{\ell-1}}\Rightarrow\node{2}{\alpha_{\ell}} &
D_4^{(3)} &\node{1}{\alpha_0}-\node{2}{\alpha_1}\Lleftarrow\node{1}{\alpha_2} \\
A_{2\ell-1}^{(2)}  &\node{1}{\alpha_1}-\node{\ver{1}{\alpha_0}}{\alpha_2}\!\!{}^2-\node{2}{\alpha_3}-\cdots-\node{2}{\alpha_{\ell-1}}\Leftarrow\node{1}{\alpha_\ell} &
E_6^{(2)} &\node{1}{\alpha_0}-\node{2}{\alpha_1}-\node{3}{\alpha_2}\Leftarrow\node{2}{\alpha_3}-\node{1}{\alpha_4}
\end{array}
\]
\caption{Twisted affine Dynkin diagrams of A,D,E.}
\label{twisted}
\end{figure}

\Def[{\cite[\S2.2]{BKM}}]
We define $\KSS{n}{s}=\frac{v^{ns}+(-1)^{s}v^{-ns}}{v^{s}+(-1)^{s}v^{-s}}(\in\MA)$ for $n,s\geq 1$ with $n\not\in 2\Z$. 
Note that $\KSS{n}{s}|_{v=1}=n$ if $s$ is odd, but $\KSS{n}{s}|_{v=1}=1$ if $s$ is even.
\edf

Recall Theorem \ref{MainTheorem} stating that for $p\geq 2$, we have 
$\SH_{\Lambda_0,w\Lambda_0-d\delta}(A^{(1)}_{p-1}) = 
\prod_{s=1}^{d}[p]_s^{N_{d,s,p}}$
for all $w\in W(A^{(1)}_{p-1})$ and $d\geq 0$
where $N_{p,d,s}=\sum_{\lambda\in\PAR(d)}\frac{m_s(\lambda)}{p-1}\prod_{u\geq 1}\binom{m_u(\lambda)+p-2}{m_u(\lambda)}$.

\Con\label{ConA22}
Let $p\geq 3$ be an odd integer and put $\ell=(p-1)/2$. We have $\SH_{\Lambda_\ell,w\Lambda_\ell-d\delta}(A^{(2)}_{p-1}) = 
\prod_{s=1}^{d}\KSSS{p}{s}{N_{\ell,d,s}}$
for all $w\in W(A^{(2)}_{p-1})$ and $d\geq 0$.
\encon

\Rem\label{Drinfeldtwisted}
Conjecture ~\ref{ConA22} should be proven similarly as Theorem ~\ref{MainTheorem} if we have
\bna
\item a proof of Drinfeld new realization of $U_v(A^{(2)}_{p-1})$ (see ~\cite{Dri}),
\item a proof of vertex operator construction of $U_v(A^{(2)}_{p-1})$-module $V(\Lambda_{(p-1)/2})$ (see ~\cite{Jin,JM}),
\item a good choice of Drinfeld loop-like generators which affords us the Lusztig lattice in terms of 
Drinfeld generators (a candidate is the PBW generators in ~\cite{BeNa}).
\ee
However, (a) is not achieved so far in general (thus, logically (b) nor (c) neither).
See also the introduction of ~\cite{Her}.
\enrem

\Rem
In virtue of ~\cite[Appendix B]{Aka} and ~\cite{Jin,JM}, we can prove Conjecture \ref{ConA22} when $p=3$
by a similar argument of ~\cite[Theorem 3]{CJ} and Theorem ~\ref{MainTheorem}.
\enrem

\Con\label{conjshap}
For a twisted affine A,D,E diagram $X=X^{(r)}_{N}$ (where $r\geq 2$), we assign the quantities $n,k\in \NNN$ and $\alpha,\beta\in\MA$ as follows.
\begin{center}
\begin{tabular}{c|ccccc}
$X^{(r)}_{N}$ & $A^{(2)}_{2n-1} (n\geq 3)$ & $A^{(2)}_{2n} (n\geq 1)$   & $D^{(2)}_{n+1} (n\geq 1)$ & $E^{(2)}_{6}$ & $D^{(3)}_{4}$ \\ \hline
$n$      & $n$        & $n$     & $n$    & $4$         & $2$         \\
$k$         & $n-1$        & $n$           & $1$             & $2$         & $1$         \\
$\alpha$    & $[2]_{n}$     & $[2n+1]^{\SUPER}$ & $[2]_{n}$      & $\{3\}_2$   & $[3]^\SUPER$ \\
$\beta$     & $[n]$       & $[2n+1]$      & $[2]$            & $[3]$       & $[2]$
\end{tabular}
\end{center}
Here for an odd integer $p\geq 1$, we define $[p]^\SUPER=(v^{p}+v^{-p})/(v+v^{-1})$.
Note that for an even integer $m\geq 2$, we have $\{p\}_m = [p]^{\SUPER}|_{v=v^{m}}$.
Put
\begin{align*}
\gamma_{X,s}=
\begin{cases}
\alpha|_{v=v^{s}} & (s\in r\Z) \\
\beta|_{v=v^{s}}  & (s\not\in r\Z),
\end{cases}\quad
f_{X,s}=
\begin{cases}
n & (s\in r\Z) \\
k    & (s\not\in r\Z),
\end{cases}
\end{align*}
then, for all $w\in W(\HA)$ and $d\geq 0$ we have
$\SH_{\Lambda_\varepsilon,w\Lambda_\varepsilon-d\delta}(\HA)=\prod_{s=1}^{d}\gamma_{X,s}^{N_{X,d,s}}$ where
\begin{align*}
N_{X,d,s} = \sum_{\lambda\in\PAR(d)}\frac{m_s(\lambda)}{f_{X,s}}\prod_{i\geq 1}\binom{f_{X,i}-1+m_{i}(\lambda)}{m_{i}(\lambda)}.
\end{align*}
\encon

\Rem
Assume $X^{(r)}_N=(B_{ij})_{i,j\in \HI}\ne A^{(2)}_{2n}$ for simplicity.
Let $X_N=(A_{i,j})_{i,j\in J}$ with $|J|=N$ and let $\mu:X_N\isoto X_N$ be a Dynkin diagram automorphism of order $r$ 
such that the folding procedure to $(X_N,\mu)$ gives $X^{(r)}_N$ (See ~\cite[\S7.9,\S8]{Kac}).
In the procedure, $I$ is identified with the set $\{\{\mu^k(j)\mid k\in \Z/r\Z\}\mid j\in J\}$ of $\mu$-orbits of $J$.
We take a map $s:I\to J$ such that $i$ corresponds to $\{\mu^k(s(i))\mid k\in \Z/r\Z\}$ for all $i\in I$.

For $i\in I$, we put $d_i = a_i^\vee/a_i$ where $a_i^\vee$ be the numerical label of the Langlands dual of $X^{(r)}_N$.
Note that we have $d_i=r/|\{\mu^k(s(i))\mid k\in \Z/r\Z\}|$ for $i\in I$.
For $t\geq 1$, we define $I(t)=\{i\in I\mid t\in d_i\Z\}$. 
Then, the quantities in the above table are designed to satisfy (see also ~\cite[\S5]{BK4})
%$\alpha$ and $\beta$ in the above table are designed to satisfy (see also ~\cite[\S5]{BK4})
$f_{X,t}=|I(t)|$ and
$\det Y^{(t)}=\gamma_{X,t}$ where $Y^{(t)}=(Y^{(t)}_{ij})_{i,j\in I(t)}$ and
\begin{align*}
Y^{(t)}_{ij}=\frac{\sum_{k\in \Z/r\Z}\exp(\frac{2\pi\sqrt{-1}}{r}tk)[A_{s(i),\mu^k(s(j))}]_t}{d_i}=
\begin{cases}
[B_{ij}]_{t} & (i=j) \\
B_{ij} & (i\ne j).
\end{cases}
\end{align*}
\enrem

%\Rem
%Conjecture \ref{conjshap} is true at least for all $n\leq 4$ and $d\leq 3$ by direct computer calculation.
%When $X^{(r)}_N=A^{(2)}_{2n},D^{(2)}_{n+1}$, Conjecture \ref{conjshap} is true at least for all $n\leq 10$ and $d\leq 10$ 
%by computer calculation using a Fock space realization of $V(\Lambda_\varepsilon)$ due to Kashiwara-Miwa-Petersen-Yung~\cite{KMPY} 
%(See also ~\cite{KKw,LT} for interpretations Kashiwara-Miwa-Petersen-Yung's construction in a more combinatorial flavour).
%\enrem

\Rem
Let $\cor$ be a field of characteristic $p\geq 2$ and let $X=X^{(r)}_{N}$ be a twisted affine A,D,E diagram.
Under the assumption that $r\ne 0$ in $\cor$,
DeConcini-Kac-Kazhdan proved that modulo-$p$ reduction of $\Gg(X)$-module $V(\Lambda_\varepsilon)$ remains irreducible
if and only if $\det X\ne 0$ in $\cor$~\cite[Remark 3.1]{DcKK}.
Conjecture \ref{conjshap} predicts an answer for Kashiwara's problem (see in \S\ref{KaPro}) 
when $\Gg=\Gg(X)$ and $\lambda(c)=1$
and it can be seen as a quantum analog of DeConcini-Kac-Kazhdan's result.
\enrem

\section{Graded Cartan determinants of the symmetric groups}\label{GRC}
Recall that in this paper, we work in the field $\cor=\Q(v)$ and its subring $\MA=\Z[v,v^{-1}]$.

\subsection{Graded representations}
\Def\label{CartanInv} 
Let $A$ be a finite dimensional algebra over a field $\corr$.
\bna
\item We denote by $\MOD{A}$ the abelian category consists of finite dimensional left $A$-modules and $A$-homomorphisms between them.
\item We define the Cartan matrix $C_A$ of $A$ to be the matrix
$([\PROC(D):D'])_{D,D'\in\IRR(\MOD{A})}\in\MAT_{|\IRR(\MOD{A})|}(\Z)$ where $\PROC(D)$ is a projective cover of $D$. 
\ee
\edf

\Def\label{GCartanInv} 
Let $A$ be a graded finite dimensional algebra over a field $\corr$, i.e., $A$ has a decomposition $A=\bigoplus_{i\in\Z}A_i$ 
into $\corr$-vector spaces such that $A_iA_j\subseteq A_{i+j}$ for all $i,j\in\Z$.
\bna
\item We denote by $\GMOD{A}$ the abelian category consists of finite dimensional left graded $A$-modules and 
degree preserving $A$-homomorphisms between them.
\item We denote by $\GPPP{A}$ the abelian category consists of finite dimensional left graded projective $A$-modules and 
degree preserving $A$-homomorphisms between them.
\item For an object $M=\bigoplus_{i\in\Z}M_i$ in $\GMOD{A}$ and $k\in\Z$, we denote by $M\langle k\rangle$ the object in $\GMOD{A}$
such that $(M\langle k\rangle)_n=M_{k+n}$ for all $n\in\Z$.
\item For $k\in\Z$, the assignment $(M\stackrel{f}{\to}N)\MAPSTO(M\langle k\rangle\stackrel{f}{\to}N\langle k\rangle)$ 
is obviously an autoequivalence both on $\GMOD{A}$ and $\GPPP{A}$ which we denote by $\langle k\rangle$.%:\GMOD{A}\iso\GMOD{A}$.
\item We define the graded Cartan pairing $\omega_A:\KKK(\GPPP{A})\times \KKK(\GMOD{A})\to\MA$ by
\begin{align*}
\textstyle\langle[P],[M]\rangle\MAPSTO \sum_{k\in\Z}\dim_{\corr}\Hom_{A}(P,M\langle k\rangle)v^k.
\end{align*}
\item We define the graded Cartan matrix $C^v_A$ of $A$ to be the matrix
\begin{align*}
\left({\textstyle\sum_{k\in\Z}[\PROC(D):D'\langle -k\rangle]v^{k}}\right)_{D,D'\in\IRR(\GMOD{A})/\SIM}\in\MAT_{|\IRR(\GMOD{A})/\SIM|}(\MA)
\end{align*}
where $\PROC(D)$ is a projective cover of $D$ in $\GMOD{A}$ and
$M\SIM N$ if and only if there exists $k\in\Z$ such that $M\langle k\rangle\cong N$ in $\GMOD{A}$.
\ee
\edf

\Rem
If $\corr$ is a splitting field for $A$, then 
$C^v_A=\TRANS{\left(\omega_A(\PROC(D),\PROC(D'))\right)_{D,D'\in\IRR(\GMOD{A})/\SIM}}$.
\enrem

\Rem
Both $\KKK(\GMOD{A})$ and $\KKK(\GPPP{A})$ have 
$\MA$-module structures %on both $\KKK(\GMOD{A})$ and $\KKK(\GPPP{A})$ 
via $v=[\langle -1\rangle]$.
$\omega_A$ is $\BAR$-sesquilinear in that $\omega_A(va,b)=v^{-1}\omega_A(a,b)$ and $\omega_A(a,vb)=v\omega_A(a,b)$ for 
all $a\in \KKK(\GPPP{A})$ and $b\in\KKK(\GMOD{A})$.
\enrem

\Rem
Let $C^v_1$ and $C^v_2$ be the graded Cartan matrices of $A$ according to
different choice of representatives of $\IRR(\GMOD{A})/\SIM$.
Then, there exists a diagonal matrix $D$ all of whose diagonal entries belong to $v^{\Z}$ such that
$C^v_2=\BAR(\TRANS{D})C^v_1D$. Thus, the graded Cartan determinant $\det C^v_A$ is well-defined.
\enrem

\Rem
It is a routine to prove that the natural map which forgets the grading gives a bijection
$\IRR(\GMOD{A})/\SIM\iso \IRR(\MOD{A})$ and $C^v_A|_{v=1}=C_{A}$.
\enrem

\subsection{The Khovanov-Lauda-Rouquier algebras}
\Def[{\cite[\S3.2.1]{Rou}}]\label{defKLR}
Let $\corr$ be a field and let $I$ be a finite set.
Take a matrix $Q=(Q_{ij}(u,v))\in \MAT_I(\corr[u,v])$ %of two variables polynomials $\cor[u,v]$ 
subject to $Q_{ii}=0,
Q_{ij}(u,v)=Q_{ji}(v,u)$
for all $i,j\in I$.
\bna
\item The KLR algebra $\KLR_n(\corr;Q)$ for $n\geq 0$ 
is an $\corr$-algebra generated by $\{x_p, \tau_a,e_{\nu}\mid 1\leq p\leq n,1\leq a<n,\nu\in I^n\}$ 
with the following defining relations for all $\mu,\nu\in I^n,1\leq p,q\leq n$, $1\leq b<a\leq n-1$.
%\vspace*{-0.32cm}
\begin{align*}\begin{tabular}{ll}
$\bullet$ $e_{\mu}e_{\nu}=\delta_{\mu\nu}e_{\mu},1=\sum_{\mu\in I^n}e_{\mu}$, $x_px_q=x_qx_p$, $x_pe_{\mu}=e_{\mu}x_p$, & 
$\bullet$ $\tau_a\tau_b=\tau_b\tau_a$ if $|a-b|>1$, \\
$\bullet$ $\tau^2_ae_{\nu}=Q_{\nu_a,\nu_{a+1}}(x_a,x_{a+1})e_{\nu}$, $\tau_ae_{\mu}=e_{s_a(\mu)}\tau_a$, & $\bullet$ $\tau_ax_p=x_p\tau_a$ if $p\ne a,a+1$, \\
$\bullet$ $(\tau_ax_{a+1}-x_a\tau_a)e_{\nu}=(x_{a+1}\tau_a-\tau_ax_a)e_{\mu}=\delta_{\nu_a,\nu_{a+1}}e_{\nu}$, & \\
\multicolumn{2}{l}{$\bullet$ $(\tau_{b+1}\tau_b\tau_{b+1}-\tau_b\tau_{b+1}\tau_b)e_{\nu}=
\delta_{\nu_b,\nu_{b+2}}((x_{b+2}-x_b)^{-1}(Q_{\nu_b,\nu_{b+1}}(x_{b+2},x_{b+1})-Q_{\nu_b,\nu_{b+1}}(x_{b},x_{b+1})))e_{\nu}$.}
\end{tabular}\end{align*}
\item For $\beta=\sum_{i\in I}\beta_i\cdot i\in\N[I]$ with $n=\HT(\beta):=\sum_{i\in I}\beta_i$, 
we define $\KLR_\beta(\corr;Q)=\KLR_n(\corr;Q)e_{\beta}$ where $e_{\beta}=\sum_{\nu\in\SEQ(\beta)}e_{\nu}$ and 
$\SEQ(\beta)=\{(i_j)_{j=1}^{n}\in I^n\mid \sum_{j=1}^{n}i_j=\beta\}$.
\item For $\lambda=\sum_{i\in I}\lambda_i\cdot i\in\N[I]$ and $\beta\in\N[I]$ with $n=\HT(\beta)$, 
we define 
\begin{align*}
\KLR^{\lambda}_n(\corr;Q) &= \KLR_n(\corr;Q)/\KLR_n(\corr;Q)(\textstyle\sum_{\nu\in I^n}x_1^{\lambda_{h_{\nu_1}}}e_{\nu})\KLR_n(\corr;Q),\\
\KLR^{\lambda}_\beta(\corr;Q) &= \KLR_\beta(\corr;Q)/\KLR_\beta(\corr;Q)(\textstyle\sum_{\nu\in \SEQ(\beta)}x_1^{\lambda_{h_{\nu_1}}}e_{\nu})\KLR_\beta(\corr;Q).
\end{align*}
\ee
\edf

As a consequence of PBW theorem for the KLR algebras~\cite[Theorem 3.7]{Rou}, 
we see that $\{e_{\beta}\mid\HT(\beta)=n\}$ 
exhausts all the primitive
central idempotents of $\KLR_n(\corr;Q)$. 
It is not difficult to see that both $\KLR^{\lambda}_n(\corr;Q)$ and $\KLR^{\lambda}_\beta(\corr;Q)$ are finite dimensional $\corr$-algebras.

\Def[{\cite[\S3.2.3]{Rou}}]
Let $A=(a_{ij})_{i,j\in I}$ be a symmetrizable {\GCM} with the symmetrization $d=(d_i)_{i\in I}$. 
Take $Q^A=(Q^A_{ij}(u,v))\in\MAT_I(\corr[u,v])$ 
subject to 
\begin{align*}
Q^A_{ii}(u,v)=0,\quad
Q^A_{ij}(u,v)=Q^A_{ji}(v,u),\quad
t_{i,j,-a_{ij},0}=t_{j,i,0,-a_{ij}}\ne 0
\end{align*}
for all $i,j\in I$ where $Q^A_{ij}(u,v)=\sum_{\substack{p,q\geq 0 \\ pd_i+qd_j=-d_ia_{ij}}} t_{ijpq}u^pv^q$.
\edf

%We abbreviate $\KLR_n(\cor;Q^A)$, $\KLR_\beta(\cor;Q^A)$, $\KLR^\lambda_n(\cor;Q^A)$, $\KLR^\lambda_\beta(\cor;Q^A)$ 
%to $\KLR_n(\cor;A)$, $\KLR_\beta(\cor;A)$, $\KLR^\lambda_n(\cor;A)$, $\KLR^\lambda_\beta(\cor;A)$ respectively.
%This abbreviation is justiced by the fact that if the Dynkin diagram of $A$ has no cycle (in the case when $A$ is of finite type
%or of affine type, this means that $A\ne A^{(1)}_{a}$ for $a\geq 1$), then the $\cor$-algebra isomorphism class
%of $\KLR_n(\cor;Q^A)$, $\KLR_\beta(\cor;Q^A)$, $\KLR^\lambda_n(\cor;Q^A)$, $\KLR^\lambda_\beta(\cor;Q^A)$ are uniquely determined.

For $n\geq 0$ and $\lambda,\beta\in\N[I]$ with $\HT(\beta)=n$, all of
$\KLR_n(\corr;Q^A)$, $\KLR_\beta(\corr;Q^A)$, $\KLR^\lambda_n(\corr;Q^A)$, $\KLR^\lambda_\beta(\corr;Q^A)$ are $\Z$-graded 
via the assignment where $\nu\in I^n,1\leq p\leq n,1\leq a<n$.
\begin{align*}
\DEG (e_{\nu})=0,\quad
\DEG (x_pe_{\nu})=2d_{\nu_p},\quad
\DEG (\tau_ae_{\nu})=-d_{\nu_a}a_{\nu_a,\nu_{a+1}}.
\end{align*}
%where $\nu\in I^n,1\leq p\leq n,1\leq a<n$.

\Th[{\cite{KKa}}]\label{KLR} Let $\corr$ be a field and let $A$ be a symmetrizable {\GCM}.
Then, we can put a $U^{\MA}_v(A)$-module structure on $\bigoplus_{n\geq 0}\KKK(\GPPP{\KLR^\lambda_n(\corr;Q^A)})$
and it is isomorphic to $V(\lambda)^{\MA}$.
\enth

\subsection{Graded Cartan determinants of the symmetric groups}\label{irrsym}
Recall that $\IRR(\MOD{\Q\mathfrak{S}_n})$ is identified with $\PAR(n)$.
For $\lambda\in\PAR(n)$ we denote by $\chi_{\lambda}$ the corresponding character.
It is characterized by $p_{\mu}=\sum_{\lambda\in\PAR(n)}\chi_{\lambda}(C_{\mu})s_{\lambda}$ in the ring of symmetric 
functions $\Lambda=\bigoplus_{n\geq 0}\varprojlim_{m\geq 0}\Z[x_1,\cdots,x_m]^{\mathfrak{S}_m}_n$.
Recall also that any field is a splitting field for $\mathfrak{S}_n$.

\Def\label{BLOCKELL}
Let $\ell\geq2$ be an integer.
\bna
\item For each $n\geq 0$, we denote by $\BLOCK_{\ell}(n)$ the set of tuple $(\rho,d)$ where $\rho$ is an $\ell$-core and $d\geq 0$ is an integer 
such that $|\rho|+\ell d=n$.
\item
For $(\rho,d)\in\BLOCK_{\ell}(n)$, we denote by $B_{\rho,d}(\subseteq\PAR(n))$ the set of partition $\lambda$ of $n$ whose
$\ell$-core of $\lambda$ is $\rho$.
\ee
\edf

The following is well-known classification of $p$-blocks of the symmetric groups due to Nakayama-Brauer-Robinson.

\Th\label{NakayamaBrauerRobinson}
Let $p\geq 2$ be a prime and let $\lambda,\mu\in\PAR(n)$ be partitons of $n\geq 0$.
The ordinary characters $\chi_{\lambda}$ and $\chi_{\mu}$ belong to the same $p$-block if and only if
$\lambda$ and $\mu$ have the same $p$-core.
\enth

\Def
Let $p\geq 2$ be a prime.
For $(\rho,d)\in \BLOCK_{p}(n)$, we define
\bna
\item $e_{\rho,d} = \sum_{\lambda\in B_{\rho,d}}\frac{\chi_{\lambda}(1)}{n!}\sum_{g\in\mathfrak{S}_{n}}\chi_{\lambda}(g)g^{-1}$,
\item $\beta_{\rho,d} = \sum_{x\in\rho}\RES(x)+d\cdot\sum_{y\in \Z/p\Z}y\in \N[\Z/p\Z]$
where $\sum_{x\in\rho}$ means that $x$ runs all the boxes of $\rho$ and $\RES(x)=j-i+p\Z$ when $x$ is located 
at $i$-th row and $j$-th column.
\ee
\edf

\Rem
$e_{\rho,d}$ is defined as an element of $\Q\mathfrak{S}_n$.
Thanks to Theorem \ref{NakayamaBrauerRobinson},
we can regard it as an element of $\mathbb{F}_p\mathfrak{S}_n$
and it is a primitive central idempotent of $\mathbb{F}_p\mathfrak{S}_n$.
\enrem

\Def\label{choiceQ}
For $\ell\geq 2$, we define $Q_{\ell}\in\MAT_{\Z/\ell\Z}(\Z[u,v])$ by
\begin{align}
(Q_{\ell})_{i,j} = 
\begin{cases}
-(u-v)^2 & \textrm{($\ell=2$ and $i\ne j$)} \\
\pm(v-u)         & \textrm{($\ell\geq 3$ and $j=i\pm 1$)} \\
1                 & \textrm{($\ell\geq 3$ and $j\ne i,i\pm1$)} \\
0                 & \textrm{(otherwise)}.
\end{cases}
\label{choiceQ}
\end{align}
For a field $\corr$, we denote by $Q_{\ell}^{\corr}\in\MAT_{\Z/\ell\Z}(\corr[u,v])$ the natural image of $Q_{\ell}$.
\edf

\Th[{\cite{Rou,BK2}}]\label{BKR} Let $\corr$ be a field of characteristic $p>0$.
Then, as $\corr$-algebras we have $e_{\rho,d}\corr\mathfrak{S}_ne_{\rho,d}\cong \KLR^{\Lambda_0}_{\beta_{\rho,d}}(\corr;Q_{p}^{\corr})$ for $(\rho,d)\in\BLOCK_{p}(n)$.
Especially, we have $\corr\mathfrak{S}_n\cong \KLR^{\Lambda_0}_n(\corr;Q_{p}^{\corr})$.
\enth

It is known that the grading comes from Theorem \ref{BKR} quantizes Ariki's 
categorification $\bigoplus_{n\geq 0}\KKK(\GP{\corr \mathfrak{S}_n})\cong V(\Lambda_0)^{\Z}$~\cite{BK1}.
Moreover, we have a compatibility of the Shapovalov form and the graded Cartan pairing.

\Th[{\cite[Theorem 4.18]{BK3}}]\label{compatibilityBK}
%Let $\cor$ be a field of characteristic $p>0$.
Let $\corr$ be an algebraically closed field and let $\ell\geq 2$.
For any $\lambda\in \MP^+$ of $A^{(1)}_{\ell-1}$, 
there exist $U^{\MA}_v(A^{(1)}_{\ell-1})$-module isomorphisms $\delta$ and $\varepsilon$ which makes
\begin{align*}
\xymatrix{
V(\lambda)^{\MA} \ar@{->}[r]^{\!\!\!\!\!\!\!\!\!\!\!\!\!\!\!\!\!\!\!\!\!\!\!\!\!\sim}_{\!\!\!\!\!\!\!\!\!\!\!\!\!\!\!\!\!\!\!\!\!\!\!\!\!\delta}\ar@{->}[d]_{a} &
\oplus_{n\geq 0}\KKK(\GPPP{\KLR^\lambda_n}) \ar@{->}[d]^{b} \\
V(\lambda)^{\MA,*} & \oplus_{n\geq 0}\KKK(\GMOD{\KLR^\lambda_n}) \ar@{->}[l]^{\!\!\!\!\!\!\!\!\!\!\!\!\!\!\!\!\!\!\!\!\!\varepsilon}_{\!\!\!\!\!\!\!\!\!\!\!\!\!\!\!\!\!\!\!\!\!\sim}}
\end{align*}
a commutative diagram. Moreover, the following diagram
\begin{align*}
\xymatrix{
V(\lambda)^{\MA}\otimes_{\MA} V(\lambda)^{\MA,*}\ar@{->}[r]^{\quad\quad\quad\quad\RSH|_{V(\lambda)^{\MA}\times V(\lambda)^{\MA,*}}}\ar@{->}_{\delta\otimes\varepsilon^{-1}}[d] & \MA \ar@{=}[d] \\
\oplus_{n\geq 0}\KKK(\GPPP{\KLR^\lambda_n})\otimes_{\MA} \oplus_{n\geq 0}\KKK(\GMOD{\KLR^\lambda_n}) \ar@{->}[r]_{\quad\quad\quad\quad\quad\quad\quad\quad\quad\oplus\omega_{\KLR^\lambda_n}} & \MA}
\end{align*}
commutes. Here we abbreviate $\KLR^\lambda_n(\corr;Q_{\ell}^{\corr})$ to $\KLR^\lambda_n$ and
\bnum
%\item We abbreviate $\KLR^\lambda_n(\corr;Q_{\ell}^{\corr})$ to $\KLR^\lambda_n$.
\item $V(\lambda)^{\MA,*}=\{v\in V(\lambda)|\forall w\in V(\lambda)^\MA, \langle v,w\rangle_{X}\in\MA\}$ where $X\in\{\SH,\QSH,\RSH\}$.
The right hand side does not depend on the choice of $X\in\{\SH,\QSH,\RSH\}$.
\item $a$ is the canonical inclusion $V(\lambda)^{\MA}\hookrightarrow V(\lambda)^{\MA,*}$.
\item $b$ is a natural map induced from 
the forgetful functor $\GPPP{\KLR^\lambda_n}\to\GMOD{\KLR^\lambda_n}$.
\ee
\enth

\Rem
The presentation of Theorem \ref{compatibilityBK} is misleading.
A crucial point of the proof is that $\KLR^\lambda_n$ is 
isomorphic to the Ariki-Koike algebra or its degeneration~\cite{BK2,Rou}.
%More precisely, for $\beta\in\N[\Z/\ell\Z]$ with $\HT(\beta)=n$, if $\KLR^\lambda_\beta\ne 0$ then
%$\KLR^\lambda_\beta$ isomorphic to a block of the Ariki-Koike algebra or that of its degeneration.
\enrem

\Def\label{DonkinC}
Let $\ell\geq 2$. We define $C^v_{\ell,d}=\QSHM_{\Lambda_0,w\Lambda_0-d\delta}(A^{(1)}_{\ell-1})$ for $w\in W(A^{(1)}_{\ell-1})$ and $d\geq 0$.
$C^{v}_{\ell,d}$ does not depend on the choice of $w$ as in Proposition \ref{Weyl} and
$C^{v}_{\ell,d}$ is related with $\RSHM_{\Lambda_0,w\Lambda_0-d\delta}(A^{(1)}_{\ell-1})$ in the sense of Proposition \ref{SHANDQSH} and Corollary \ref{detichi}.
\edf

By combining Proposition \ref{SHANDQSH}, Corollary \ref{detichi}, Theorem \ref{MainTheorem} and Theorem \ref{compatibilityBK},
we get the graded Cartan determinants of the symmetric groups.

\Th\label{gradeddeter2}
Let $\corr$ be a filed of characteristic $p>0$.
For each $(\rho,d)\in \BLOCK_{p}(n)$, we have
$\det C^v_{e_{\rho,d}\corr\mathfrak{S}_ne_{\rho,d}}=\det C^v_{\ell,d}=\prod_{s=1}^{d}[p]_s^{N_{p,d,s}}$ 
where
\begin{align*}
N_{p,d,s}=\sum_{\lambda\in\PAR(d)}\frac{m_s(\lambda)}{p-1}\prod_{u\geq 1}\binom{m_u(\lambda)+p-2}{m_u(\lambda)}
=\sum_{(\lambda_i)_{i=1}^{p-1}\in\PAR_{p-1}(d)}m_s(\lambda_1).
\end{align*}
\enth

\subsection{Graded Cartan determinants of the Iwahori-Hecke algebras of type A}\label{IwahoriHecke}
\Def
Let $\MAA$ be a commutative domain and take $q\in\MAA$.
The Iwahori-Hecke algebra of type A with a quantum parameter $q$
is an $\MAA$-algebra $\MH_n(R;q)$ generated by $\{T_i\mid 1\leq i<n\}$ and obeys the following defining relations.
\begin{align*}
(T_b+1)(T_b-q)=0,\quad T_aT_{a+1}T_a=T_{a+1}T_aT_{a+1},\quad T_bT_c=T_cT_b
\end{align*}
where $1\leq a\leq n-2$ and $1\leq b,c<n$ with $|b-c|>1$.
\edf

\Def
For $\ell\geq 2$, we fix a field $\corrr{\ell}$ which has an $\ell$-th primitive root of unity $\qq{\ell}$.
\edf

\Rem\label{blockiw}
For $\ell\geq 2$ and $n\geq 0$, it is known that 
\bnum
\item $\corrr{\ell}$ is a splitting field for $\MH_n(\corrr{\ell};\qq{\ell})$ (see ~\cite[\S2.2]{Don}),
\item\label{blockiw2} $\BLOCK_{\ell}(n)$ 
parametrizes the set of all primitive central idempotents of $\MH_n(\corrr{\ell};\qq{\ell})$.
\ee
\enrem

Concerning Remark \ref{blockiw} (\ref{blockiw2}), for $(\rho,d)\in \BLOCK_{\ell}(n)$ we denote by $e'_{\rho,d}$ the
corresponding primitive central idempotent of $\MH_n(\corrr{\ell};\qq{\ell})$ as in ~\cite[\S5]{DJ}.
The choice is the same as in ~\cite{Rou,BK2} based on the Jucy-Murphy elements.

\Th[{\cite{Rou,BK2}}]\label{BKR2} Let $\ell\geq 2$.
Then, as $\corrr{\ell}$-algebras we 
have $e'_{\rho,d}\MH_n(\corrr{\ell};\qq{\ell})e'_{\rho,d}\cong \KLR^{\Lambda_0}_{\beta_{\rho,d}}(\corrr{\ell};Q_{\ell}^{\corrr{\ell}})$ for $(\rho,d)\in\BLOCK_{p}(n)$.
Especially, we have $\MH_n(\corrr{\ell};\qq{\ell})\cong \KLR^{\Lambda_0}_n(\corrr{\ell};Q_{\ell}^{\corrr{\ell}})$.
\enth

By the same arguments as for Theorem \ref{gradeddeter2}, we get the following.

\Th\label{gradeddeter3}
For $\ell\geq 2$ and $(\rho,d)\in\BLOCK_{\ell}(n)$,
we have $\det C^v_{e'_{\rho,d}\MH_n(\corrr{\ell};\qq{\ell})e'_{\rho,d}}%=\det C^v_{\ell,d}
=\prod_{s=1}^{d}[\ell]_s^{N_{\ell,d,s}}$.
\enth

\Rem\label{qAriki}
Let $\ell\geq 2$ and assume $\CHAR\corrr{\ell}=0$.
In virtue of the graded version of Ariki's theorem~\cite[Corollary 5.15]{BK3} and the graded Brauer-Humphreys 
reciprocity~\cite[Theorem 2.17]{HM}, we have
$C^v_{e'_{\rho,d}\MH_n(\corrr{\ell};\qq{\ell})e'_{\rho,d}}=\TRANS{D_{\rho,d}}D_{\rho,d}$
where $D_{\rho,d}=(d_{\lambda,\mu}(v))_{{\lambda\in B_{\rho,d},\mu\in B_{\rho,d}\cap \RPAR_{\ell}(d)}}$ and $\RPAR_{\ell}(d)$ is the set
of $\ell$-regular partitions of $d$.
Recall that we have the Fock space representation $\mathcal{F}=\bigoplus_{\lambda\in\PAR}\cor |\lambda\rangle$ of $U_v(A^{(1)}_{\ell-1})$ due to Hayashi
and $\bigsqcup_{n\geq 0}\{G(\mu):=\sum_{\lambda\in\PAR(n)}d_{\lambda,\mu}(v)|\lambda\rangle\mid \mu\in\RPAR_{\ell}(n)\}$ is 
Lusztig's canonical basis (also known as Kashiwara's lower global basis) 
of $V(\Lambda_0)\cong U_v(A^{(1)}_{\ell-1})|\emptypartition\rangle$ (see ~\cite[\S6]{LLT}).
\enrem

\section{Conjectures on the graded Cartan invarinats of the symmetric groups}\label{GRI}
\subsection{Generalized blocks of the symmetric groups}
In their paper ~\cite{KOR},
K\"ulshammer-Olsson-Robinson initiated the study of generalized blocks for finite groups. % in ~\cite{KOR}.

\Def[{\cite[\S1]{KOR}}]\label{combblock}
Let $G$ be a finite group and let $\MCC\subseteq G$ be a subset of $G$ invariant under
conjugation, i.e., a union of conjugacy classes of $G$. %which is invariant under conjugation.
Let $\KORC(G)$ be the set of complex-valued class functions on $G$ with scalar product $\langle\alpha,\beta\rangle_{\MCC}=\frac{1}{|G|}
\sum_{g\in\MCC}\alpha(g)\overline{\beta(g)}$.
\bna
\item For two ordinary characters $\psi,\varphi$ of $G$, we say that $\psi$ and $\varphi$ are directly $\MCC$-linked (and 
denote by $\psi\approx_\MCC\varphi$) if
$\langle\psi,\varphi\rangle_{\MCC}\ne 0$. Clearly, the relation $\approx_\MCC$ is symmetric.
\item A $\MCC$-block of $G$ is an equivalence class of the ordinary characters of $G$ under the equivalence
relation $\simeq_{\MCC}$ which is defined as the transitive closure of the relation $\approx_{\MCC}$.
\ena
\edf

\Rem
For $G_{p'}:=\{g\in G\mid \ORD_G(g)\not\in p\Z\}$ where $\ORD_G(g)$ is the order of $g$ in $G$, the notion of $G_{p'}$-blocks 
coincides with the usual notion of $p$-blocks.
\enrem

\Def
Let $\ell\geq 2$ be an integer.
\bna
\item A partition $\lambda$ is called $\ell$-class regular if no parts of $\lambda$ is 
divisible by $\ell$ (i.e., $m_{k\ell}(\lambda)=0$ for all $k\geq 1$). We denote by
$\CPAR_{\ell}(n)$ the set of $\ell$-class regular partitions of $n$.
\item An $\ell$-block of $\mathfrak{S}_n$ is defined as a $\MCC_{\ell'}$-block 
for $\MCC_{\ell'}=\bigsqcup_{\lambda\in \CPAR_{\ell}(n)}C_{\lambda}(\subseteq \mathfrak{S}_n)$ where $C_{\lambda}$ is a conjugacy class of $\mathfrak{S}_n$ 
consists of elements whose cycle type is $\lambda$.
\ena
\edf

\Th[{\cite[Theorem 5.13]{KOR}}]\label{KORmainthm}
Let $\ell\geq 2$ be an integer and let $\lambda,\mu\in\PAR(n)$ be partitons of $n\geq 0$.
The ordinary characters $\chi_{\lambda}$ and $\chi_{\mu}$ (see \S\ref{irrsym}) belong to the same $\ell$-block if and only if
$\lambda$ and $\mu$ have the same $\ell$-core.
\enth

\subsection{Conjectures of K\"ulshammer-Olsson-Robinson and Hill}
For each $\MCC$-block $B$ of $G$ (see Definition \ref{combblock}) where $\MCC$ is closed (i.e.,
$\langle x\rangle=\langle y\rangle$ implies $y\in \MCC$ for any $x\in \MCC$ and $y\in G$), K\"ulshammer-Olsson-Robinson assigned the Cartan group.
%``generalized Cartan matrix'' which 
%is a $|B|\times |B|$ integer matrix up to unimodular equivalence over $\Z$~\cite[\S1]{KOR}.

\Def[{\cite[\S1]{KOR}}]\label{CARTANGROUP} Let $G$ be a finite group and let $\MCC$ %$=\sqcup_{i=1}^{r}C_i$ 
be a union of %$r$ 
conjugacy classes of $G$ which is closed.
%such that $\langle x\rangle=\langle y\rangle$ implies $y\in \MCC$ for any $x\in \MCC$ and $y\in G$. 
\bna
\item We define two additive subgroups $\KORR_G(\MCC)$ and $\KORP_G(\MCC)$ of $\KORC(G)$ by
\begin{align*}
\KORR_G(\MCC) &= \SPAN_\Z\{\chi^{\MCC}\mid \chi\in\IRR(\C G)\},\\
\KORP_G(\MCC) &= \KORR_G(\MCC)\cap \SPAN_\Z\{\chi\mid \chi\in\IRR(\C G)\}
\end{align*}
where $\IRR(\C G)$ is the set of ordinary characters of $G$ and
$f^{\MCC}\in \KORC(G)$ is the class function that takes the same value as $f$
on $\MCC$ and vanishes elsewhere for $f\in \KORC(G)$.
\item The Cartan group $\CART_G(\MCC)$ is $\KORR_G(\MCC)/\KORP_G(\MCC)$ which is a finite group.
\ee
\edf

We call the invariant factors of $\CART_G(\MCC)$ the generalized Cartan invariants.
When $\ell=p$ is a prime and $\MCC=G_{p'}$, the generalized Cartan invariants 
are the usual Cartan invariants (i.e., the elementary divisors of the Cartan matrix $C_{\overline{\F_p}G}$). 
Recall that in the $p$-modular representation theory, Cartan invariants carry an information of certain
centralizer groups.

\Th[{\cite[Part III,\S16]{BrNe}}]\label{classicalmodular}
Let $G$ be a finite group and let $p\geq2$ be a prime.
Take a set of representatives $g_1,\cdots,g_k$ of the $p$-regular conjugacy classes $C_1,\cdots C_k$ of $G$.
The Cartan invariants of $\overline{\F_p}G$ is given by the multiset $\{|\SYL_p(C_G(g_i))|\mid 1\leq i\leq k\}$.
\enth

In the case of $G=\mathfrak{S}_n$ and $\MCC=\MCC_{\ell'}$ for $\ell\geq 2$,
we may identify the Cartan group as $\COKER(\Z^{\oplus|\RPAR_{\ell}(n)|}\to\Z^{\oplus|\RPAR_{\ell}(n)|},\XXX\mapsto\XXX C_{\MH_n(\corrr{\ell};\qq{\ell})})$ because 
of the following.

\Th[{\cite[\S2.2]{Don}}]
Let $\ell\geq 2$. For any $n\geq 0$, there exist isomorphisms $\theta_n$ and $\zeta_n$ such that
the following diagram commutes.
\begin{align*}
\xymatrix{
\KKK(\PROJ{\MH_n(\corrr{\ell};\qq{\ell})})\ar@{->}[r]^{\quad\quad\sim}_{\quad\quad\theta_n}\ar@{->}[d] & \KORP_{\mathfrak{S}_n}(\MCC_{\ell'})\ar@{^{(}->}[d] \\
\KKK(\MOD{\MH_n(\corrr{\ell};\qq{\ell})})\ar@{->}[r]^{\quad\quad\sim}_{\quad\quad\zeta_n} & \KORR_{\mathfrak{S}_n}(\MCC_{\ell'})}
\end{align*}
%Moreover, 
%for $x\in \KKK(\PROJ{\MH_n(\corrr{\ell};\qq{\ell})})$ and $y\in\KKK(\MOD{\MH_n(\corrr{\ell};\qq{\ell})})$ we have
%\begin{align*}
%\omega_{\MH_n(\corrr{\ell};\qq{\ell})}(x,y)=\langle\theta_n(x),\zeta_n(y)\rangle_{\MCC_{\ell'}}.
%\end{align*}
\enth

\Def
Let $\MAA$ be a commutative ring with $1$.
\bna
\item We say that two $m\times m$ matrices $X,Y\in\MAT_m(\MAA)$ are unimodular equivalent over $\MAA$ (and denote by $X\CONG{\MAA} Y$) 
if there exist $P,Q\in \GL_m(\MAA)$ such that $X=PYQ$.
\item For a matrix $X$ and a multiset $S$ both  consist of elements of $\MAA$,
we abbreviate $X\CONG{\MAA} \DIAG(S)$ to $X\CONG{\MAA} S$
where $\DIAG(S)=(s_{ij})_{i,j\in I}$ is a 
diagonal matrix with a multiset identity $\{s_{ii}\mid i\in I\}=S$.
\ee
\edf

\Def
Let $n\geq 1$ and let $\Pi$ be a subset of the set of all prime numbers.
\bna
\item We define $\PRS(n)$ to be the set of prime divisors of $n$ (when $n=1$, $\PRS(n)=\emptyset$).
\item We define $n_{\Pi}$ to be the $\Pi$-part of $n$, i.e.,
the unique positive integer $n_{\Pi}$ such that $n\in n_{\Pi}\Z$ and $\PRS(n_{\Pi})\subseteq \Pi, \PRS(n/n_{\Pi})\cap \Pi=\emptyset$.
\ee
\edf

For $\ell\geq 2$ and $(\rho,d)\in\BLOCK_{\ell}(n)$, we put $C_{\ell,d}:=C^v_{\ell,d}|_{v=1}$ (see Definition \ref{DonkinC}).

\Con[{\cite[Conjecture 6.4]{KOR}}]\label{KORcon}
Let $\ell\geq2$ be an integer. We define
\begin{align*}
r_{\ell}(\lambda)=\prod_{k\in\NNN\setminus \ell\Z}({\ell}/{\GCD(\ell,k)})^{\lfloor \frac{m_k(\lambda)}{\ell} \rfloor}\cdot 
\bigl\lfloor \frac{m_k(\lambda)}{\ell} \bigr\rfloor!_{\PRS({\ell}/{\GCD(\ell,k)})}
\end{align*}
for a partition $\lambda$.
Then, for $n\geq 0$, we have %(see also Remark \ref{depth})
\begin{align*}
\bigoplus_{(\rho,d)\in\BLOCK_{\ell}(n)}C_{\ell,d} \CONG{\Z}
\{ r_{\ell}(\lambda)\mid \lambda\in\CPAR_{\ell}(n) \}.
\end{align*}
\encon

\Rem\label{depth}
For $n\geq 0$, we may regard $\bigoplus_{(\rho,d)\in\BLOCK_{\ell}(n)}C^v_{\ell,d}$ as
\begin{align*}
\GRM_{\cor,\BAR,\MA}\left(\oplus_{\substack{\mu\in P(\lambda) \\ \HT(\Lambda_0-\mu)=n}}V(\lambda)_{\mu},
\oplus_{\substack{\mu\in P(\lambda) \\ \HT(\Lambda_0-\mu)=n}}\QSH,V(\lambda)^{\MA}_{\mu}\right)=
\oplus_{\substack{\mu\in P(\lambda) \\ \HT(\Lambda_0-\mu)=n}}\QSHM_{\Lambda_0,\mu}(A^{(1)}_{\ell-1}).
\end{align*}
\enrem

\Rem
Let $p\geq 2$ be a prime. Then,
\bna
\item $p$-regular conjugacy classes of $\mathfrak{S}_n$ is given by $\{C_{\lambda}\mid \lambda\in\CPAR_p(n)\}$,
\item the centralizer group $C_{\mathfrak{S}_n}(g)$ for $g\in C_{\lambda}$ is given by $\bigoplus_{k\geq 1}(\Z/k\Z)\wr \mathfrak{S}_{m_k(\lambda)}$.
\ee
\enrem

Let $\ell=\prod_{p\in\PRS(\ell)}p^r$ be a prime factorization of an integer $\ell\geq 2$.
Hill settled affirmatively the KOR conjecture when $r\leq p$ for each $p\in\PRS(\ell)$~\cite[Theorem 1.3]{Hil}.
In the course of the proof, Hill proposed a following refinent of the KOR conjecture into each $\ell$-blocks for $|\PRS(\ell)|=1$.

\Def
Let $n\geq 1$ and let $\ell\geq 2$.
We define $a_{\ell}(n)\geq 1$ and $\nu_{\ell}(n)\geq 0$
to be the unique integers such that $a_{\ell}(n)\ell^{\nu_{\ell}(n)}=n$ and $a_{\ell}(n)\not\in\ell\Z$.
\edf

\Con[{\cite[Conjecture 10.5]{Hil}}]\label{Hcon}
Let $p\geq 2$ be a prime and let $r\geq 1$ be an integer. 
For a partition $\lambda$ we define a power of $p$ by
\begin{align*}
\log_p I_{p,r}(\lambda) = \sum_{n\in\NNN\setminus p^r\Z}((r-\nu_p(n))m_n(\lambda)+\sum_{t\geq 1}\lfloor m_n(\lambda)/p^t \rfloor).
\end{align*}
Put $\ell=p^r$. Then, for each $d\geq 0$ 
we have 
\begin{align*}
C_{\ell,d} \equiv_{\Z} \bigsqcup_{s=1}^{d}\bigsqcup_{\lambda\in\PAR(s)}\{I_{p,r}(\lambda)\}^{u(\ell-2,d-s)}.
\end{align*}
\encon

Hill proved that his conjecture implies the KOR conjecture.
In fact, they are equivalent. % (see Corollary \ref{equiva}).

\Th[{\cite[Theorem 1.2]{Hil}}]
Let $\ell\geq 2$ be an integer and define $I_{\ell}(\lambda)=\prod_{p\in\PRS(\ell)}I_{p,\nu_p(\ell)}(\lambda)$ for a partition $\lambda$.
If Conjecture \ref{Hcon} is true, then for each $d\geq 0$ we have
\begin{align*}
C_{\ell,d} \equiv_{\Z} \bigsqcup_{s=1}^{d}\bigsqcup_{\lambda\in\PAR(s)}\{I_{\ell}(\lambda)\}^{u(\ell-2,d-s)}.
\end{align*}
\enth

\Cor\label{equiva}
Conjecture \ref{KORcon} and Conjecture \ref{Hcon} are equivalent.
\encor

\Proof
Let $\ell\geq 2$ be an integer and let us assume Conjecture \ref{Hcon} is true.
Then, we see that Conjecture \ref{KORcon} is true since
Bessenrodt-Hill proved the multiset identity
\begin{align}
\{ r_{\ell}(\lambda)\mid \lambda\in\CPAR_{\ell}(n) \}=
\bigsqcup_{(\rho,d)\in\BLOCK_{\ell}(n)}\bigsqcup_{s=1}^{d}\bigsqcup_{\lambda\in\PAR(s)}\{I_{\ell}(\lambda)\}^{u(\ell-2,d-s)}
\label{BHmulti}
\end{align}
for each $n\geq 0$~\cite[Theorem 5.2]{BH}. 

Conversely, assume Conjecture \ref{KORcon} is true. For a prime $p\geq 2$ and an integer $r\geq 1$, we put $\ell=p^r$.
It is enough to prove by induction on $w\geq 0$ thatwe have
$C_{\ell,w} \equiv_{\Z} \bigsqcup_{s=1}^{w}\bigsqcup_{\lambda\in\PAR(s)}\{I_{p,r}(\lambda)\}^{u(\ell-2,w-s)}$.
When $w=0$, it is trivial. Because we have 
\begin{align*}
\{ r_{\ell}(\lambda)\mid \lambda\in\CPAR_{\ell}(\ell w) \}
\CONG{\Z}
\bigoplus_{(\rho,d)\in\BLOCK_{\ell}(\ell w)}C_{\ell,d}
=C_{\ell,w}\oplus
\bigoplus_{\substack{(\rho',d)\in\BLOCK_{\ell}(\ell w) \\ d<w}}C_{\ell,d},
\end{align*}
inductively we must have $C_{\ell,w} \equiv_{\Z} \bigsqcup_{s=1}^{w}\bigsqcup_{\lambda\in\PAR(s)}\{I_{p,r}(\lambda)\}^{u(\ell-2,w-s)}$ in virtue of (\ref{BHmulti}).
Here we use an elementary linear algebra fact that, for a PID $\MAA$ and finitely generated $\MAA$-modules $X,Y$ and $Y'$,
$X\oplus Y\cong X\oplus Y'$ as $R$-modules implies $Y\cong Y'$.
\QED

%Thus, in order to prove KOR conjecture, it is enough to prove Hill's conjecture\footnote{Bessenrodt-Hill 
%proved that Hill's conjecture and KOR-conjecture are equivalent.}.

\subsection{Gradation of Hill's conjecture}
For integers $a\geq 0$ and $b\geq 1$, we denote by $a\% b$ the remainder of $a$ by $b$, namely
the unique integer $0\leq c<b$ such that $a-c\in b\Z$.

\Con\label{myconj}
Let $p\geq 2$ be a prime and let $r\geq 1$ be an integer. We define
\begin{align*}
I^{v}_{p,r}(\lambda) = \prod_{n\in\NNN\setminus p^r\Z}\prod_{k=1}^{m_n(\lambda)}[p^{r+\nu_p(k)-\nu_p(n)}]_{a_p(k)p^{\nu_p(n)}}
\end{align*}
for a partition $\lambda$.
Put $\ell=p^r$. Then, for each $d\geq 0$ we have 
\begin{align*}
C^{v}_{\ell,d}\CONG{\MA} \bigsqcup_{s=1}^{d}\bigsqcup_{\lambda\in\PAR(s)}\{I^{v}_{p,r}(\lambda)\}^{u(\ell-2,d-s)}.
\end{align*}
\encon

\Rem
We have $I^{v}_{p,r}(\lambda)|_{v=1}=I_{p,r}(\lambda)$ for any $\lambda\in\PAR$.
Thus, Conjecture \ref{myconj} implies Conjecture \ref{Hcon}.
Note that we do not necessarily have $I^{v}_{p,r}(\lambda_1)\in \MA I^{v}_{p,r}(\lambda_2)$ or $I^{v}_{p,r}(\lambda_2)\in \MA I^{v}_{p,r}(\lambda_1)$
for $\lambda_1,\lambda_2\in\PAR(n)$ unlike at $v=1$.
\enrem

\Rem\label{computer_check}
In the setting of Conjecture \ref{myconj},
Conjecture \ref{myconj} implies the following. % weak version %of Conjecture \ref{myconj}
\begin{align}
C^{v}_{\ell,d}\CONG{\Q[v,v^{-1}]} \bigsqcup_{s=1}^{d}\bigsqcup_{\lambda\in\PAR(s)}\{I^{v}_{p,r}(\lambda)\}^{u(\ell-2,d-s)}.
\label{weakversion}
\end{align}
When $r=1$,
(\ref{weakversion}) is the same as ~\cite[Conjecture 8.2 (i)]{ASY} (see Remark \ref{qAriki}).
Although we could not prove (\ref{weakversion}), it seems tractable by an elementary method
since we have
$C^{v}_{\ell,d}\CONG{\Q[v,v^{-1}]} \bigsqcup_{s=1}^{d}\bigsqcup_{\lambda\in\PAR(s)}\{{\textstyle\prod_{i\geq 1}[\ell]_i^{m_i(\lambda)}}\}^{u(\ell-2,d-s)}$
which follows from the proof of Proposition \ref{comp}.
\enrem

As a support of Conjecture \ref{myconj},
we prove that Conjecture \ref{myconj} gives the graded Cartan determinants of Theorem \ref{MainTheorem}.

\Th\label{conjcheck}
Let $p\geq 2$ be a prime and let $r\geq 1$ be an integer.
Put $\ell=p^r$. Then, for each $d\geq 0$ we have the equality 
\begin{align*}
\prod_{s=1}^{d}[\ell]_{s}^{N_{\ell,d,s}} = \prod_{s=1}^{d}\prod_{\lambda\in\PAR(s)}I^{v}_{p,r}(\lambda)^{u(\ell-2,d-s)}.
\end{align*}
\enth

\Rem
Ando-Suzuki-Yamada proved in ~\cite[Theorem 4.3]{ASY} that for any $\ell\geq 2$
\begin{align}
\prod_{s=1}^{d}[\ell]_{s}^{N_{\ell,d,s}} = \prod_{s=1}^{d}\prod_{\lambda\in\PAR(s)}Q_{\ell}(\lambda)^{u(\ell-2,d-s)}.
\label{ASYformula}
\end{align}
where
$Q_{\ell}(\lambda):=\prod_{n\in\NNN\setminus\ell\Z}\prod_{k=1}^{m_{n}(\lambda)}[\ell^{1+\nu_\ell(k)}]_{a_{\ell}(k)}$.
Note that $I^{v}_{p,r}(\lambda)=Q_{p^r}(\lambda)$ when $r=1$. Thus, Theorem \ref{conjcheck} gives another
product expansion of the graded Cartan determinants of the symmetric groups conjecturally comes from 
nice representatives of the unimodular equivalence classes of the graded Cartan matrices (but only when $|\PRS(\ell)|=1$).
\enrem

\Proof
In virtue of (\ref{ASYformula}), it is enough to prove 
\begin{align}
\prod_{s=1}^{d}\prod_{\lambda\in\PAR(s)}Q_{p^r}(\lambda)^{u(\ell-2,d-s)}
=
\prod_{s=1}^{d}\prod_{\lambda\in\PAR(s)}I^{v}_{p,r}(\lambda)^{u(\ell-2,d-s)}.
\label{mokuhyo_conjcheck}
\end{align}
Since $[a^b]_c=\prod_{i=1}^{b}[a]_{ca^{i-1}}$ for $a,b,c\geq 1$ and $\nu_{p^r}(n)=\lfloor\nu_{p}(n)/r\rfloor,a_{p^r}(n)=a_{p}(n)p^{\nu_{p}(n)\%r}$ for $n\geq 1$,
it is enough to show the following multiset identity in order to prove (\ref{mokuhyo_conjcheck})
\begin{align}
\begin{split}
{} &{\quad}
\bigsqcup_{s=1}^{d}\bigsqcup_{\lambda\in\PAR(s)}
\bigsqcup_{n\in\NNN\setminus p^r\Z}\bigsqcup_{k=1}^{m_n(\lambda)}\{a_p(k)p^t\mid \nu_p(k)\%r\leq t<r+\nu_p(k)\}^{u(\ell-2,d-s)}\\
&=
\bigsqcup_{s=1}^{d}\bigsqcup_{\lambda\in\PAR(s)}
\bigsqcup_{n\in\NNN\setminus p^r\Z}\bigsqcup_{k=1}^{m_n(\lambda)}\{a_p(k)p^t\mid \nu_p(n)\leq t<r+\nu_p(k)\}^{u(\ell-2,d-s)}.
\label{mokuhyo_conjcheck2}
\end{split}
\end{align} 
Clearly, (\ref{mokuhyo_conjcheck2}) follows from the following multiset identity
\begin{align*}
{} &{\quad}
\bigsqcup_{\lambda\in\PAR(s)}
\bigsqcup_{n\in\NNN\setminus p^r\Z}\bigsqcup_{k=1}^{m_n(\lambda)}\{a_p(k)p^t\mid \nu_p(n)\leq t<r+\nu_p(k)\} \\
&=
\bigsqcup_{\lambda\in\PAR(s)}
\bigsqcup_{n\in\NNN\setminus p^r\Z}\bigsqcup_{k=1}^{m_n(\lambda)}\{a_p(k)p^t\mid \nu_p(k)\%r\leq t<r+\nu_p(k)\}
\end{align*}
for $s\geq 1$ which follows from Proposition \ref{tsaigo}.\QED

\Prop\label{tsaigo}
Let $p\geq 2$ be a prime and let $r\geq 1$ be an integer.
We have the following multiset identity for each $d\geq 1$ and $u\in\NNN\setminus p\Z$.
\begin{align*}
\bigsqcup_{\lambda\in\PAR(d)}
\bigsqcup_{n\geq1}\bigsqcup_{\substack{1\leq k\leq m_n(\lambda) \\ a_p(k)=u}}\{\nu_p(n)\} 
=
\bigsqcup_{\lambda\in\PAR(d)}
\bigsqcup_{n\geq1}\bigsqcup_{\substack{1\leq k\leq m_n(\lambda) \\ a_p(k)=u}}\{\nu_p(k)\%r\}.
\end{align*}
\enprop

Before giving a proof, we need preparatory definitions.

\Def
For $d\geq 1$, we define maps $\CUT_{d},\INFL_{d}:\PAR\to\PAR$ by the following.
\bna
\item for $k\geq 1$, $m_k(\CUT_{d}(\lambda))=m_k(\lambda)$ if $k\not\in d\Z$, otherwise $m_k(\CUT_{d}(\lambda))=0$.
\item $\ell(\lambda)=\ell(\INFL_{d}(\lambda))$ and $(\INFL_{d}(\lambda))_i=d\lambda_i$ for $1\leq i\leq \ell(\lambda)$,
\ee
\edf

\Def\label{defofP}
Let $p\geq 2$ be a prime and let $r\geq 1$ be an integer.
Note that any integer $a\geq 1$ can be uniquely expressed as $a=\sum_{i=0}^{r-1}a_ip^i$ with $a_{r-1}\geq 0$ and $0\leq a_j<p$ for $0\leq j\leq r-2$.
Define the subset $\PPP^{p,r}_a$ of $\PAR(a)$ consists of all the partitions $\nu$ such that
\bnum
\item\label{ine} for all $0\leq k<r$, we have $\sum_{h=k}^{r-1}a_hp^{h-k}\geq \sum_{h=k}^{r-1}m_{p^h}(\nu)p^{h-k}$,
\item we have $\sum_{h=0}^{r-1}a_hp^{h}= \sum_{h=0}^{r-1}m_{p^h}(\nu)p^{h}$.
\ee
\edf

\Rem\label{CHARACTERIZATIONPPP}
$\PPP^{p,r}_a$ is characterized as the smallest subset of $\PAR(a)$ such that %that satisfies the following conditions.
\bna
\item Define $\lambda(a)\in\PAR(a)$ by $m_{p^i}(\lambda(a))=a_i$ for $0\leq i<r$. 
Then, we have $\lambda(a)\in \PPP^{p,r}_a$.
\item For any $\lambda\in \PPP^{p,r}_a$ and any $0\leq i<r$ with $m_i(\lambda)>0$, define $\mu\in\PAR(a)$ to be 
%a partition defined by
\begin{align*}
m_j(\mu) = \begin{cases}
m_i(\lambda)-1 & \textrm{($j=i$)} \\
m_{i-1}(\lambda)+p & \textrm{($j=i-1$)} \\
m_j(\lambda) & \textrm{(otherwise)}
\end{cases}
\end{align*}
for $0\leq j<r$. Then, we have $\mu\in \PPP^{p,r}_a$.
\ee

Thus, we have $\PPP^{p,r}_a=\{\lambda\in\CPAR_{p^r}(a)\mid \forall f\in \NNN\setminus p^{\N},m_f(\lambda)=0\}$.
\enrem

\Proof
Because of Lemma \ref{bunkaitoP}, it is enough to prove the following multiset identity for each $a\geq 1$ and $u\in\NNN\setminus p\Z$.
\begin{align}
\bigsqcup_{\lambda\in \PPP^{p,r}_a}
\bigsqcup_{n\geq1}\bigsqcup_{\substack{1\leq k\leq m_n(\lambda) \\ a_p(k)=u}}\{\nu_p(n)\} 
=
\bigsqcup_{\lambda\in \PPP^{p,r}_a}
\bigsqcup_{n\geq1}\bigsqcup_{\substack{1\leq k\leq m_n(\lambda) \\ a_p(k)=u}}
\{\nu_p(k)\%r\}.
\label{shimeshitaishiki}
\end{align}

For any $0\leq s<r$, we define sets $U_1(s)$ and $U_2(s)$ as follows.
\begin{align*}
U_1(s) &= \{(\lambda,j)\in \PPP^{p,r}_a\times\N \mid up^j\leq m_{p^s}(\lambda)\},\\
U_2(s) &= \{(\lambda,j,t)\in \PPP^{p,r}_a\times\N\times\{0,\cdots,r-1\} \mid up^j\leq m_{p^t}(\lambda),j-s\in r\Z\}.
\end{align*}

Note that (\ref{shimeshitaishiki}) is the same as $\bigsqcup_{0\leq s<r}\{s\}^{|U_1(s)|}=\bigsqcup_{0\leq s<r}\{s\}^{|U_2(s)|}$.
Thus, it is enough to construct a bijection between $U_1(s)$ and $U_2(s)$ for any $0\leq s<r$.

We claim that the map $U_2(s)\to U_1(s),(\lambda,j,t)\mapsto (\mu,i)$ given by
\begin{align*}
i = t+j-s,\quad
m_{p^w}(\mu)=
\begin{cases}
m_{p^t}(\lambda)-up^j & \textrm{($s\ne t$ and $w=t$)} \\
m_{p^s}(\lambda)+up^i & \textrm{($s\ne t$ and $w=s$)} \\
m_{p^w}(\lambda)      & \textrm{(otherwise)}
\end{cases}
\end{align*}
where $0\leq w<r$ is well-defined.
To see this, only non-trivial check should be to prove $\mu\in \PPP^{p,r}_a$ when $s>t$.
In order to prove this, it is enough to check (\ref{ine}) in Definition \ref{defofP} for $\nu=\mu$ and $h=s$. Since $\lambda\in \PPP^{p,r}_a$, we have
$\sum_{h=t}^{r-1}a_hp^{h-t}\geq \sum_{h=t}^{r-1}m_{p^h}(\lambda)p^{h-t}$. Thus, %we have
\begin{align*}
\sum_{h=s}^{r-1}a_hp^{h-s}+p^{t-s}\sum_{h=t}^{s-1}a_hp^{h-t} 
&\geq \sum_{h=s}^{r-1}m_{p^h}(\lambda)p^{h-s}+p^{t-s}\sum_{h=t}^{s-1}m_{p^h}(\lambda)p^{h-t} \\
&\geq \sum_{h=s}^{r-1}m_{p^h}(\lambda)p^{h-s}+m_{p^t}(\lambda)p^{h-t} \geq \sum_{h=s}^{r-1}m_{p^h}(\mu)p^{h-s}.
\end{align*}
Since $0\leq p^{t-s}\sum_{h=t}^{s-1}a_hp^{h-t}\leq p^{t-s}\sum_{h=t}^{s-1}(p-1)\cdot p^{h-t}=1-p^{t-s}<1$, we are done.

Similarly, we see that the map $U_1(s)\to U_2(s),(\mu,i)\mapsto (\lambda,j,t)$ given by
\begin{align*}
t=i\%r,\quad
j=i+s-t,\quad
m_{p^w}(\lambda)=
\begin{cases}
m_{p^t}(\mu)+up^j & \textrm{($s\ne t$ and $w=t$)} \\
m_{p^s}(\mu)-up^i & \textrm{($s\ne t$ and $w=s$)} \\
m_{p^w}(\mu)      & \textrm{(otherwise)}
\end{cases}
\end{align*}
is well-defined. Clearly, it is a converse of the map $U_2(s)\to U_1(s)$ above.
\QED

\Lemma\label{bunkaitoP}
Let $p\geq 2$ be a prime and let $r\geq 1$ be an integer.
For each $d\geq 1$, the multiset $S^{p,r}_{d}:=\{\CUT_{p^r}(\lambda)\mid\textrm{$\lambda\in\PAR(d)$ such that $\CUT_{p^r}(\lambda)\ne\emptypartition$}\}$
has a decomposition of the form
$S^{p,r}_{d}=\bigsqcup_{i\in I}\SORT(\prod_{j\in J_i}\INFL_{d_{i,j}}(\PPP^{p,r}_{a_{i,j}}))$ such that 
\bna
\item\label{hoge1} $I$ is a finite non-empty set and $\{J_i\mid i\in I\}$ is a family of finite non-empty sets,
\item\label{hoge2} $a_{i,j}\geq 1$ and $d_{i,j}\in\NNN\setminus p\Z$ for $i\in I$ and $j\in J_i$,
\item\label{hoge3} for any $i\in I$ and $j\ne j'\in J_i$, take $\lambda\in \INFL_{d_{i,j}}(\PPP^{p,r}_{a_{i,j}})$ and $\mu\in \INFL_{d_{i,j'}}(\PPP^{p,r}_{a_{i,j'}})$. 
Then, we have either $m_f(\lambda)=0$ or $m_f(\mu)=0$ for all $f\geq 1$,
\item $\SORT:\PAR^\ast\to\PAR$ is a map that maps a finite sequence of partitions $(\lambda^{(g)})_{g\in G}$ %as its argument
to a partition $\mu$ such that $\{\mu_f\mid 1\leq f\leq \ell(\mu)\}=\bigsqcup_{g\in G}\{\lambda^{(g)}_{f}\mid 1\leq f\leq \ell(\lambda^{(g)})\}$
as multisets.
\ee
\enlemma

\Proof
Let $\lambda$ be a non-empty partition.
Clearly, $L=\{1,\cdots,\ell(\lambda)\}$ has a decomposition of the form $L=\bigsqcup_{k\in K}L_k$ such that
\bnum
\item for any $k\in K$, we have $L_k\ne\emptyset$ and 
there exists a unique $d_k\in \NNN\setminus p\Z$ such that $\lambda_j\in d_k\Z$ and $\lambda_j/d_k\in p^\N$
for all $j\in L_k$,
\item for any $k\ne k'\in K$, we have $d_k\ne d_{k'}$.
\ee

Put $\lambda^{(k)}=\SORT(((\lambda_j/d_k))_{j\in L_k})$ for $k\in K$. 
As in Remark \ref{CHARACTERIZATIONPPP}, there exists a unique $a_k\geq 1$ such that $\lambda^{(k)}\in \PPP^{p,r}_{a_k}$.
Thus, we have $\lambda\in\SORT(\prod_{k\in K}\INFL_{d_k}(\PPP^{p,r}_{a_k}))$.

From this construction, we see that $\CPAR_{p^r}(d)$ has a decomposition of the form
$\CPAR_{p^r}(d)=\bigsqcup_{i\in I}\SORT(\prod_{j\in J_i}\INFL_{d_{i,j}}(\PPP^{p,r}_{a_{i,j}}))$ that satisfies (\ref{hoge1}), (\ref{hoge2}) and (\ref{hoge3})
for each $d\geq 1$.

Because we have
$S^{p,r}_{d}=\bigsqcup_{0\leq e< d/p^r}\CPAR_{p^r}(d-p^re)^{|\PAR(e)|}$, we are done.
\QED

\subsection{Gradation of K\"ulshammer-Olsson-Robinson's conjecture when $\ell=p^r$}
By the same reason to consider a gradation of Hill's conjecture (see \S\ref{introinv}),
we want to consider a correct gradation of the KOR conjecture.
Though we could not find it in general, we propose a gradation of the KOR conjecture only when $|\PRS(\ell)|=1$.

\Con\label{myconj2}
Let $p\geq 2$ be a prime and let $r\geq 1$ be an integer. Put $\ell=p^r$.
Then, for each $n\geq 0$, we have 
$\bigoplus_{(\rho,d)\in\BLOCK_{\ell}(n)}C^{v}_{\ell,d} \CONG{\MA}
\{ r^{v}_{p,r}(\lambda)\mid \lambda\in\CPAR_{\ell}(n) \}$
where %we define $r^{v}_{p,r}(\lambda)$ for any partition $\lambda$ as follows.
\begin{align*}
r^{v}_{p,r}(\lambda)
=\prod_{k\geq 1}\prod_{t=1}^{\lfloor\frac{m_k(\lambda)}{\ell}\rfloor}[p^{r-\nu_p(k)}]_{tp^{\nu_p(k)}}\cdot[p^{\nu_p(t)}]_{a_p(t)p^{\nu_p(k)}}
=\prod_{k\geq 1}\prod_{t=1}^{\lfloor\frac{m_k(\lambda)}{\ell}\rfloor}[p^{r-\nu_p(k)+\nu_p(t)}]_{a_p(t)p^{\nu_p(k)}}.
\end{align*}
\encon

\Rem
When $r=1$, Conjecture \ref{myconj2} implies ~\cite[Conjecture 8.2 (ii)]{ASY} which predicts 
$\bigoplus_{(\rho,d)\in\BLOCK_{p}(n)}C^{v}_{p,d} \CONG{\Q[v,v^{-1}]}
\{ r^{v}_{p,1}(\lambda)\mid \lambda\in\CPAR_{p}(n) \}$ for $n\geq 0$ and a prime $p\geq 2$.
\enrem

\Prop\label{conjequiv}
Conjecture \ref{myconj} implies Conjecture \ref{myconj2}. % are equivalent.
More precisely, we have the following multiset identity for each $r,n\geq 1$ and a prime $p\geq 2$ where $\ell=p^r$.
\begin{align*}
\bigsqcup_{(\rho,d)\in\BLOCK_{\ell}(n)}\bigsqcup_{s=1}^{d}\bigsqcup_{\lambda\in\PAR(s)}\{I^{v}_{p,r}(\lambda)\}^{u(\ell-2,d-s)}
=
\{ r^{v}_{p,r}(\lambda)\mid \lambda\in\CPAR_{\ell}(n) \}.
\end{align*}
\enprop

\Proof
Since $r^{v}_{p,r}(\lambda)=I^{v}_{p,r}(\RED_{\ell}(\lambda))$ and $I^{v}_{p,r}(\lambda)=I^{v}_{p,r}(\CUT_{\ell}(\lambda))$ for any partition $\lambda$ where
$m_k(\RED_{\ell}(\lambda))=\lfloor m_k(\lambda)/\ell\rfloor$ for all $k\geq 1$,
Proposition \ref{conjequiv} follows from Lemma \ref{saigo2}.
%You can see a proof of Lemma \ref{saigo2} in the proof of ~\cite[Theorem 3.3]{ASY}\footnote{Their statement ``$\#\alpha^{-1}(\nu)=\#\beta^{-1}(\nu)$''
%in the proof of ~\cite[Theorem 3.3]{ASY} is the same thing as Lemma \ref{saigo2}.}.
\QED

\Lemma[{\cite[Lemma 5.5]{BH}}]\label{saigo2}
%Let $p\geq 2$ be a prime and let $r\geq 1$ be an integer.
%Put $\ell=p^r$. We have the following multiset identity for each $d\geq 1$.
Let $\ell\geq 2$ be an integer. We have the multiset identity %for any $n\geq 0$.
\begin{align*}
\bigsqcup_{(\rho,d)\in\BLOCK_{\ell}(n)}\bigsqcup_{s=1}^{d}\bigsqcup_{\lambda\in\PAR(s)}\{\CUT_{\ell}(\lambda)\}^{u(\ell-2,d-s)}
=
\{ \RED_{\ell}(\lambda)\mid \lambda\in\CPAR_{\ell}(n) \}
\end{align*}
for any $n\geq 0$.
\enlemma

\end{document}